\documentclass[11pt]{article}
\usepackage{amssymb}
\usepackage{times}

\title{$\omega$-powers and descriptive set theory.\indent}
\author{Dominique LECOMTE}
\date{\it ~J. Symbolic Logic~\rm 70, 4 (2005), 1210-1232}

\newcommand{\Ana}{{\it\Sigma}^{1}_{1}}
\newcommand{\Ca}{{\it\Pi}^{1}_{1}}
\newcommand{\Borel}{{\it\Delta}^{1}_{1}}

\newcommand{\Bormone}{{\it\Pi}^{0}_{1}}
\newcommand{\Bormtwo}{{\it\Pi}^{0}_{2}}
\newcommand{\Boraone}{{\it\Sigma}^{0}_{1}}
\newcommand{\ana}{{\bf\Sigma}^{1}_{1}}
\newcommand{\ca}{{\bf\Pi}^{1}_{1}}
\newcommand{\borel}{{\bf\Delta}^{1}_{1}}
\newcommand{\boraone}{{\bf\Sigma}^{0}_{1}}
\newcommand{\boratwo}{{\bf\Sigma}^{0}_{2}}

\newcommand{\boraxi}{{\bf\Sigma}^{0}_{\xi}}

\newcommand{\bortwo}{{\bf\Delta}^{0}_{2}}
\newcommand{\borone}{{\bf\Delta}^{0}_{1}}

\newcommand{\bormone}{{\bf\Pi}^{0}_{1}}
\newcommand{\bormtwo}{{\bf\Pi}^{0}_{2}}

\newcommand{\bormxi}{{\bf\Pi}^{0}_{\xi}}

\newcommand{\bormm}{{\bf\Pi}^{0}_{m}}

\newtheorem{thm} {Theorem}
\newtheorem{defi} [thm] {Definition}
\newtheorem{cor} [thm] {Corollary}
\newtheorem{lem} [thm] {Lemma}

\newtheorem{prop} [thm] {Proposition}

\begin{document}

\maketitle

\noindent {\footnotesize {\bf Abstract.} We study the sets of the infinite sentences 
constructible with a dictionary over a finite alphabet, from the 
viewpoint of descriptive set theory. Among other things, this gives some 
true co-analytic sets. The case where the dictionary is finite is 
studied and gives a natural example of a set at the level $\omega$ of the 
Wadge hierarchy.}\bigskip\smallskip

\noindent\bf {\Large 1 Introduction.}\rm\bigskip

 We consider the finite alphabet $n=\{0,\ldots ,n-1\}$, where $n\geq 
2$ is an integer, and a dictionary over this alphabet, i.e., a subset $A$ of the 
set $n^{<\omega}$ of finite words with letters in $n$.

\begin{defi} The $\omega\mbox{-}power$ associated to $A$ is the set $A^\infty$ of the infinite sentences 
constructible with $A$ by concatenation. So we have 
${A^\infty := \{{a_{0}}{a_{1}}\ldots\in n^\omega /\forall i\!\in\!\omega~~a_{i}\!\in\! A\}}$.\end{defi} 

 The $\omega$-powers play a crucial role in the 
characterization of subsets of $n^\omega$ accepted by finite automata 
(see Theorem 2.2 in [St1]). We will study these objects from the 
viewpoint of descriptive set theory. The reader should see [K1] for the classical results of this theory; we will also use the notation of this book. The questions we study are the following:\bigskip

\noindent (1) What are the possible levels of topological complexity for 
the $\omega$-powers? This question was asked by P. Simonnet in [S], and studied in [St2]. O. Finkel (in [F1]) and A. Louveau 
proved independently that $\ana$-complete $\omega$-powers exist.  O. Finkel 
proved in [F2] the existence of a $\bormm$-complete $\omega$-power for 
each integer $m\geq 1$.\bigskip

\noindent (2) What is the topological complexity of the set of 
dictionaries whose associated $\omega$-power is of a given level of 
complexity? This question arises naturally when we look at the characterizations of $\bormone$, $\bormtwo$ and $\boraone$ $\omega$-powers obtained in [St2] (see Corollary 14 and 
Lemmas 25, 26).\bigskip

\noindent (3) We will recall that an $\omega$-power is an analytic subset 
of $n^\omega$. What is the topological complexity of the set of codes 
for analytic sets which are $\omega$-powers? This question was asked 
by A. Louveau. This question also makes sense for the set of codes 
for $\boraxi$ (resp., $\bormxi$) sets which are $\omega$-powers. 
And also for the set of codes for Borel sets which are 
$\omega$-powers.\bigskip

 As usual with descriptive set theory, the point is not only the 
computation of topological complexities, but also the hope that these 
computations will lead to a better understanding of the studied 
objects. Many sets in this paper won't be clopen, in particular won't 
be recursive. This gives undecidability results.

\vfill\eject

\noindent  $\bullet$ We give the answer to Question (2) for the very first levels 
($\{\emptyset\}$, its dual class and $\borone$). This contains a 
study of the case where the dictionary is finite. In particular, we 
show that the set of dictionaries whose associated $\omega$-power is 
generated by a dictionary with two words is a 
$\check D_{\omega}(\boraone)$-complete set. This is a surprising result because this complexity 
is not clear at all on the definition of the set.\bigskip

\noindent $\bullet$ We give two proofs of the fact that the relation 
``$\alpha\in A^\infty$" is $\ana$-complete. One of these proofs is used later 
to give a partial answer to Question (2). To understand this answer, the reader should see 
[M] for the basic notions of effective descriptive set theory. Roughly 
speaking, a set is effectively Borel (resp., effectively Borel in $A$) if its 
construction based on basic clopen sets can be coded with a recursive (resp., 
recursive in $A$) sequence of integers. This answer is the\bigskip

\noindent\bf Theorem.\it\ The following sets are true 
co-analytic sets:\medskip

\noindent - $\{A\in 2^{n^{<\omega}}/A^\infty\in\Borel (A)\}$.\smallskip

\noindent - $\{A\in 2^{n^{<\omega}}/A^\infty\in\boraxi\cap\Borel (A)\}$, for 
$1\leq\xi <\omega_{1}$.\smallskip

\noindent - $\{A\in 2^{n^{<\omega}}/A^\infty\in\bormxi\cap\Borel (A)\}$, for 
$2\leq\xi <\omega_{1}$.\rm\bigskip

 This result also comes from an analysis of Borel 
$\omega$-powers: $A^\infty$ is Borel if and only if we can choose in a 
Borel way the decomposition of any sentence of $A^\infty$ into words of 
$A$ (see Lemma 13). This analysis is also related to Question (3) and to some 
Borel uniformization result for $G_{\delta}$ sets locally with Borel 
projections. We will specify these relations.\bigskip

\noindent  $\bullet$ A natural ordinal rank can be defined on the 
complement of any $\omega$-power, and we study it; its knowledge 
gives an upper bound of the complexity of the $\omega$-power.\bigskip

\noindent  $\bullet$ We study the link between Question (1) and the extension 
ordering on finite sequences of integers.\bigskip

\noindent $\bullet$ Finally, we give some examples of $\omega$-powers complete 
for the classes $\borone$, $\boraone\oplus\bormone$, 
$D_{2}(\boraone )$, $\check D_{2}(\boraone )$, $\check D_{3}(\boraone )$ and 
$\check D_{2}(\boratwo )$.\bigskip\smallskip

\noindent\bf {\Large 2 Finitely generated $\omega$-powers.}\rm\bigskip

\bf\noindent Notation.\rm ~In order to answer to Question (2), we set 
$${\bf\Sigma}_{0} := \{A\subseteq n^{<\omega}/A^\infty 
=\emptyset\}\mbox{,}~~ 
{\bf\Pi}_{0} := \{A\subseteq n^{<\omega}/A^\infty =n^\omega\}\mbox{,}$$
$${\bf\Delta}_{1} := \{A\subseteq 
n^{<\omega}/A^\infty\in\borone \}\mbox{,}$$
$${\bf\Sigma}_{\xi} := \{A\subseteq n^{<\omega}/A^\infty\in\boraxi 
\}\mbox{,}~~{\bf\Pi}_{\xi} := \{A\subseteq n^{<\omega}/A^\infty\in\bormxi 
\}~~~(\xi\geq 1)\mbox{,}$$
$${\bf\Delta} := \{A\subseteq n^{<\omega}/A^\infty\in\borel \}.$$
$\bullet$ If $A\subseteq n^{<\omega}$, then we set $A^{-}:=A\setminus\{\emptyset\}$.\bigskip

\noindent $\bullet$ We define, for $s\in n^{<\omega}$ and $\alpha\in n^\omega$,
$\alpha -s := \big(\alpha (\vert s\vert ),\alpha (\vert s\vert +1),...\big)$.\bigskip

\noindent $\bullet$ If ${\cal S}\subseteq (n^{<\omega})^{<\omega}$, then we set 
${\cal S}^{*}:=\{S^{*}:=S(0)\ldots S(\vert s\vert -1)/S\in {\cal S}\}$.

\vfill\eject 

\noindent $\bullet$ We define a recursive map $\pi : n^\omega\times\omega^\omega\times
\omega\rightarrow n^{<\omega}$ by 
$$\pi (\alpha ,\beta ,q):=\left\{\!\!\!\!\!\!\!\!
\begin{array}{ll} 
& {\big(\alpha\big(0\big),\ldots,\alpha (\beta [0])\big)}~\mbox{if}~q=0\mbox{,}\cr 
& \big(\alpha (1+\Sigma_{j<q}~\beta [j]),...,\alpha (\Sigma_{j\leq q}~
\beta [j])\big)~\mbox{otherwise.}
\end{array}
\right.$$
 We always have the following equivalence:
$$\alpha\!\in\! A^\infty ~\Leftrightarrow 
~~\exists\beta\!\in\!\omega^\omega ~~[\big(\forall m\! >\! 0~~\beta (m)\! >\! 0\big)~\mbox{and}~\big(\forall q\!\in\!\omega~~\pi (\alpha ,\beta ,q)\!\in\! A\big)].$$

\begin{prop}~([S]) $A^\infty\in\ana$ for all 
$A\subseteq n^{<\omega}$. If $A$ is finite, then $A^\infty\in\bormone$.\end{prop} 

\bf\noindent Proof.\rm ~We define a continuous map $c:(A^-)^\omega\rightarrow 
n^\omega$ by the formula $c\big((a_{i})\big):={a_{0}}{a_{1}}\ldots$ We have 
$A^\infty = c[(A^{-})^\omega]$, and $(A^-)^\omega$ is a Polish space (compact if 
$A$ is finite).$\hfill\square$

\begin{prop} If $A^\infty\in\borone$, then there exists a 
finite subset $B$ of $A$ such that $A^\infty = B^\infty$.\end{prop} 

\bf\noindent Proof.\rm ~Set $E_{k}:= 
\{\alpha\in n^\omega /\alpha\lceil 
k\in A~\mbox{and}~\alpha -\alpha\lceil k\in A^\infty\}$. It is an open subset of  
$n^\omega$ since $A^\infty$ is open, and 
$A^\infty\subseteq \bigcup_{k>0} E_{k}$. We can find an integer $p$ such that 
$A^\infty\subseteq \bigcup_{0<k\leq p} E_{k}$, by compactness of $A^\infty$. 
Let $B:=A\cap n^{\leq p}$. If $\alpha\in A^\infty$, then we can find an 
integer $0<k_{0}\leq p$ such that $\alpha\lceil k_{0}\in A~{\rm 
and}~\alpha -\alpha\lceil k_{0}\in A^\infty$. Thus $\alpha\lceil 
k_{0}\in B$. Then we do it again with $\alpha -\alpha\lceil k_{0}$, 
and so on. Thus we have $\alpha\in B^\infty =A^\infty$.$\hfill\square$\bigskip

\bf\noindent Remark.\rm ~This is not true if we only assume that $A^\infty$ 
is closed. Indeed, we have the following counter-example, due to O. Finkel:
$$A:=\{s\in 2^{<\omega}/\forall i\leq\vert s\vert ~~2.\mbox{Card}(\{j<i/s(j)=1\})
\geq i\}.$$
We have $A^\infty = \{\alpha\in 2^\omega /\forall i\in\omega ~~2.\mbox{Card}(\{j<i/\alpha (j)=1\})\geq i\}$ and if $B$ is finite and $B^\infty 
=A^\infty$, $B\subseteq A$ and $101^20^2\ldots\notin B^\infty$.

\begin{thm} (a) ${\bf\Sigma}_{0}=\{\emptyset,\{\emptyset\}\}$ is 
$\bormone$-complete.\smallskip

\noindent (b) ${\bf\Pi}_{0}$ is a dense $\boraone$ subset of $2^{n^{<\omega}}$. In 
particular, ${\bf\Pi}_{0}$ is $\boraone$-complete.\smallskip

\noindent (c) ${\bf\Delta}_{1}$ is a $K_{\sigma}\setminus\bormtwo$ subset of 
$2^{n^{<\omega}}$. In particular, ${\bf\Delta}_{1}$ is $\boratwo$-complete.\end{thm}

\bf\noindent Proof.\rm ~(a) Is clear.\bigskip

\noindent (b) If we can find $m\in\omega$ with $n^m\subseteq A$, then $A^\infty=
n^\omega$. As ${\{A\subseteq n^{<\omega}/\exists m\in\omega ~n^m\subseteq A\}}$ 
is a dense open subset of $2^{n^{<\omega}}$, the density follows. The formula 
$$A\in {\bf\Pi}_{0} ~\Leftrightarrow ~\exists m~\forall s\in 
n^m~\exists q\leq m~~s\lceil q\in A^{-}$$
shows that ${\bf\Pi}_{0}$ is $\Boraone$, and comes from Proposition 3.\bigskip 

\noindent (c) If $A^\infty\in\borone$, then we can find $p>0$ such that 
${A^\infty = (A\cap n^{\leq p})^\infty}$, by Proposition 3. 
So let $s_{1}, \ldots ,s_{k},t_{1}, \ldots ,t_{l}\in 2^{<\omega}$ be 
such that $A^\infty = \bigcup_{1\leq i\leq k} N_{s_{i}} = 
{n^\omega\setminus (\bigcup_{1\leq j\leq l} N_{t_{j}})}$. For each 
$1\leq j\leq l$, and for each sequence $s\in [(A^{-})^{<\omega}]^{*}\setminus\{\emptyset\}$, 
$t_{j}\not\prec s$. So we have 
$$A^\infty\!\in\!\borone~\Leftrightarrow\left\{\!\!\!\!\!\!
\begin{array}{ll} 
& \exists p>0~~\exists k,l\!\in\!\omega~~\exists s_{1},\ldots ,s_{k},t_{1},\ldots ,
t_{l}\!\in\! 2^{<\omega}~~\bigcup_{1\leq i\leq k} N_{s_{i}} = 
{n^\omega\!\setminus\! (\bigcup_{1\leq j\leq l} N_{t_{j}})}\cr 
& \mbox{and}~~\forall 1\leq j\leq l~~\forall s\!\in\! [(A^{-})^{<\omega}]^{*}\!
\setminus\!\{\emptyset\}~~t_{j}\not\prec s~~\mbox{and}~~\forall\alpha\!\in\! n^\omega\cr 
& \{\alpha\!\notin\!\bigcup_{1\leq i\leq k}\!\! N_{s_{i}}~\mbox{or}~\exists\beta\!\in\! 
p^\omega~[\big(\forall m\! >\! 0~~\beta (m)\! >\! 0\big)~\mbox{and}~\big( \forall 
q\!\in\!\omega~~\pi (\alpha ,\beta ,q)\!\in\! A\big)]\}.
\end{array}
\right.$$
 This shows that ${\bf\Delta}_{1}$ is a $K_{\sigma}$ subset of $2^{n^{<\omega}}$.
 
\vfill\eject 
 
 To show that it is not $\bormtwo$, it is enough to see that its 
intersection with the closed set 
$${\{A\! \subseteq\!  n^{<\omega}/A^\infty\! \not=\!  n^\omega\}}$$ 
is dense and co-dense in this closed set (see (b)), by Baire's theorem. 
So let $O$ be a basic clopen subset of $2^{n^{<\omega}}$ meeting this 
closed set. We may assume that it is of the form 
$$\{A\subseteq n^{<\omega}/\forall i\leq k~~s_{i}\in A~\mbox{and}~\forall j
\leq l~~t_{j}\notin A\}\mbox{,}$$ where $s_{0},\ldots ,s_{k},t_{0}, \ldots ,t_{l}\in n^{<\omega}$ 
and $\vert s_{0}\vert >0$. Let $A:= \{s_{i}/i\leq k\}$. Then $A\in O$ and 
$A^\infty$ is in $\bormone\setminus\{\emptyset,n^\omega\}$. There are two cases.\bigskip 

 If $A^\infty\!\in\!\borone$, then we have to find $B\!\in\! O$ with 
$B^\infty\!\notin\!\borone$. Let $u_{0}, \ldots ,u_{m}\in n^{<\omega}$  
with $\bigcup_{p\leq m} N_{u_{p}} = n^\omega\setminus A^\infty$. Let 
$r\in n\setminus\{ u_{0}(\vert u_{0}\vert -1)\}$, $s:=u_{0} r^{\vert u_{0}\vert +\mbox{max}_{j\leq l}~\vert t_j\vert }$ and $B := A\cup\{s\}$. Then $B\in O$ and 
$s^\infty\in B^\infty$. Let us show that $s^\infty$ is not in the 
interior of $B^\infty$. Otherwise, we could find an integer $q$ such 
that $N_{s^q}\subseteq B^\infty$. We would have 
$\alpha := {s^q}{u_{0}}u_{0}\big(\vert u_{0}\vert -1\big)r^\infty\in B^\infty$. As 
$N_{u_{0}}\cap A^\infty = \emptyset$, the 
decomposition of $\alpha$ into nonempty words of $B$ would start with 
$q$ times $s$. If this decomposition could go on, then we would have $u_{0}=
\big( u_{0}(\vert u_{0}\vert -1)\big) ^{\vert u_{0}\vert }$. Let $v\in n^{<\omega}$ 
be such that $N_{v}\subseteq A^\infty$. We have 
$v\big( u_{0}(\vert u_{0}\vert -1)\big)^\infty\in A^\infty$, so 
$\big( u_{0}(\vert u_{0}\vert -1)\big)^\infty\in N_{u_{0}}\cap A^\infty$. 
But this is absurd. Therefore $B^\infty\notin\borone$.\bigskip
    
 If $A^\infty\!\notin\!\borone$, then we have to find $B\!\in\! O$ such that 
$B^\infty\!\in\!\borone\!\setminus\!\{n^\omega\}$. Notice that 
${n^\omega\!\not=\!\bigcup_{i\leq k} N_{s_{i}}}$. So let $v\in n^{<\omega}$ 
be non constant such that $N_{v}\cap \bigcup_{i\leq k} N_{s_{i}}=\emptyset$. 
We set 
$${D\! :=\! A\cup\bigcup_{r\in n\setminus\{ v(0)\}}~\!\{ (r)\}\!\cup\!\{ {v(0)}^
{\vert v\vert }\}}\mbox{,}$$ 
${B\! :=\! A\cup\{s\!\in\! n^{<\omega}/\vert s\vert \! >\!\mbox{max}_{j\leq l}\vert t_j\vert ~\mbox{and}~
\exists t\!\in\! D~t\prec s\}}$. We get ${B^\infty = \bigcup_{t\in D} N_{t}\in\borone}$ and 
$${N_{v}\cap B^\infty =\emptyset}\mbox{,}$$ 
so ${B^\infty\not=n^\omega}$.$\hfill\square$\bigskip

 Now we will study ${\cal F}:=\{A\subseteq n^{<\omega}/\exists B\subseteq 
n^{<\omega}~\mbox{finite}~A^\infty =B^\infty\}$. 

\begin{prop} ${\cal F}$ is a 
co-nowhere dense $\boratwo$-hard subset of $2^{n^{<\omega}}$.\end{prop} 

\bf\noindent Proof.\rm ~By Proposition 3, 
if $A^\infty=n^\omega$, then there exists an integer $p$ such that 
${A^\infty = (A\cap n^{\leq p})^\infty}$, so ${\bf\Pi}_{0}\subseteq 
{\cal F}$ and, by Theorem 4, ${\cal F}$ is co-nowhere dense. We 
define a continuous map ${\phi : 2^\omega\!\rightarrow\!  
2^{n^{<\omega}}}$ by the formula $\phi (\gamma )\! :=\!\{0^k1/\gamma (k)\! =\! 1\}$. 
If ${\gamma\!\in\! P_{f}\! :=\!\{\alpha\!\in\! 2^\omega /\exists p~\forall m\!
\geq\! p~~\alpha (m)\! =\! 0\}}$, then $\phi (\gamma )\in {\cal F}$. If 
$\gamma\notin P_{f}$, then the concatenation map is an homeomorphism from 
$\phi (\gamma )^\omega$ onto $\phi (\gamma )^\infty$, thus $\phi (\gamma )^\infty$ 
is not $K_{\sigma}$. So $\phi (\gamma )\notin {\cal F}$, by Proposition 2. 
Thus the preimage of ${\cal F}$ by $\phi$ is $P_{f}$, and ${\cal F}$ is 
$\boratwo$-hard.$\hfill\square$\bigskip

 Let ${\cal G}_{p}:=\{A\subseteq n^{<\omega}/\exists s_{1},\ldots ,s_{p}\in 
n^{<\omega}~A^\infty =\{s_{1},\ldots ,s_{p}\}^\infty\}$, so that ${\cal F}=
\bigcup_{p}~{\cal G}_{p}$. We have ${\cal G}_{0}={\bf\Sigma}_{0}$, so 
${\cal G}_{0}$ is $\Bormone\setminus\boraone$.

\vfill\eject

\begin{prop} ${\cal G}_{1}$ is $\bormone\setminus\boraone$. In particular, ${\cal G}_{1}$ is $\bormone$-complete.\end{prop}

\bf\noindent Proof.\rm ~If $p\in\omega\setminus\{ 0\}$, then 
$\{0,1^p\}\notin {\cal G}_{1}$ since $B^\infty =\{s^\infty\}$ if 
$B=\{s\}$. Thus $\{0\}$ is not an interior point of ${\cal G}_{1}$ 
since the sequence $(\{0,1^p\})_{p>0}$ tends to $\{0\}$. So 
${\cal G}_{1}\notin\boraone$.\bigskip

\noindent $\bullet$ Let $(A_{m})\subseteq {\cal G}_{1}$ tending to $A\subseteq 
n^{<\omega}$. If $A\subseteq \{\emptyset\}$, then 
$A^\infty =\emptyset = \{\emptyset\}^\infty$, so $A\in {\cal G}_{1}$. If 
$A\not\subseteq \{\emptyset\}$, then let $t\in A^{-}$ and $\alpha_{0}:=t^\infty$. There 
exists an integer $m_{0}$ such that $t\in A_{m}$ for $m\geq m_{0}$. Thus we may 
assume that $t\in A_{m}$ and $A_{m}^\infty\not=\emptyset$. So let $s_{m}\in 
n^{<\omega}\setminus\{\emptyset\}$ be such that $A_{m}^\infty =\{s_{m}\}^\infty =
\{s_{m}^\infty\}$. We have $s_{m}^\infty=\alpha_{0}$. Let $b:=\mbox{min}\{a\in\omega\setminus\{0\} /(\alpha_{0}\lceil a)^\infty=\alpha_{0}\}$.\bigskip

\noindent $\bullet$ We will show that $A_{m}\subseteq 
\{(\alpha_{0}\lceil b)^q/q\in\omega\}$. Let $s\in 
A_{m}\setminus\{\emptyset\}$. As $s^\infty =\alpha_{0}$, we can find an 
integer $a>0$ such that $s=\alpha_{0}\lceil a$, and $b\leq a$. Let 
$r<b$ and $q$ be integers so that $a=q.b+r$. We have, if $r>0$,
$$\begin{array}{ll}
\alpha_{0}\!\!\! 
& =(\alpha_{0}\lceil a)^\infty = (\alpha_{0}\lceil b)^\infty = (\alpha_{0}\lceil q.b)
(\alpha_{0}\lceil a-\alpha_{0}\lceil q.b)\alpha_{0}\cr 
& = {(\alpha_{0}\lceil b)^q}(\alpha_{0}\lceil a-\alpha_{0}\lceil q.b)\alpha_{0} =
(\alpha_{0}\lceil a-\alpha_{0}\lceil q.b)\alpha_{0} =
{(\alpha_{0}\lceil r)}\alpha_{0} = {(\alpha_{0}\lceil r)}^\infty .
\end{array}$$
 Thus, by minimality of $b$, $r=0$ and we are done.\bigskip
 
\noindent $\bullet$ Let $u\in A$. We can find an integer $m_{u}$ such that $u\in 
A_{m}$ for $m\geq m_{u}$. So there exists an integer $q_{u}$ such 
that $u=(\alpha_{0}\lceil b)^{q_{u}}$. Therefore 
$A^\infty =\{(\alpha_{0}\lceil b)^\infty\}=\{\alpha_{0}\lceil b\}^\infty$ 
and $A\in {\cal G}_{1}$.$\hfill\square$\bigskip
 
\noindent\bf Remark.\rm ~Notice that this shows that we can find $w\in 
n^{<\omega}\setminus\{\emptyset\}$ such 
that $A\subseteq\{ w^q/q\in\omega\}$ if $A\in {\cal G}_{1}$. Now we study 
${\cal G}_{2}$. The next lemma is just Corollary 6.2.5 in [Lo]. 
   
\begin{lem} Two finite sequences which commute are 
powers of the same finite sequence.\end{lem}

\bf\noindent Proof.\rm ~Let $x$ and $y$ be finite sequences 
with $xy=yx$. Then the subgroup of the free group on $n$ generators 
generated by $x$ and $y$ is abelian, hence isomorphic to $\mathbb{Z}$. One 
generator of this subgroup must be a finite sequence $u$ such that $x$ 
and $y$ are both powers of $u$.$\hfill\square$

\begin{lem} Let $A\in {\cal G}_{2}$. Then there exists 
a finite subset $F$ of $A$ such that $A^\infty=F^\infty$.\end{lem} 

\bf\noindent Proof.\rm ~We will show more. Let 
$A\notin {\cal G}_{1}$ satisfying $A^\infty=\{s_{1},s_{2}\}^\infty$, with 
$\vert s_{1}\vert \leq\vert s_{2}\vert $. Then\bigskip

\noindent (a) The decomposition of $\alpha$ into words of $\{s_{1},s_{2}\}$ is unique for each 
$\alpha\in A^\infty$ (this is a consequence of Corollaries 6.2.5 and 6.2.6 in [Lo]).\smallskip

\noindent (b) ${s_{2}}{s_{1}}\perp {s_{1}}^q{s_{2}}$ for each integer $q>0$, 
and ${s_{2}}{s_{1}}\wedge {s_{1}}^q{s_{2}}={s_{1}}{s_{2}}\wedge {s_{2}}{s_{1}}$.\smallskip

\noindent (c) $A\subseteq [\{s_{1},s_{2}\}^{<\omega}]^{*}$.\bigskip

\noindent $\bullet$ We prove the first two points. We split into cases.\bigskip

\noindent 2.1. $s_{1}\perp s_{2}$.\bigskip

 The result is clear.
 
\vfill\eject
 
\noindent 2.2. $s_{1}\prec_{\not=} s_{2}\not\prec s_{1}^\infty$.\bigskip

 Here also, the result is clear (cut $\alpha$ into words of length $\vert s_{1}\vert $).\bigskip

\noindent 2.3. $s_{1}\prec_{\not=} s_{2}\prec s_{1}^\infty$.\bigskip

 We can write $s_{2}={s_{1}^m}s$, where $m>0$ and $s\prec_{\not=}s_{1}$. Thus 
$s_{2}s_{1}={s_{1}^m}ss_{1}$ and ${s_{1}^{m+1}}s\prec s_{1}^qs_{2}$ if 
$q>0$. But ${s_{1}^m}ss_{1}\perp {s_{1}^{m+1}}s$ otherwise 
$ss_{1}=s_{1}s$, and $s$, $s_{1}$ $s_{2}$ would be powers of some 
sequence, which contradicts $A\notin {\cal G}_{1}$.\bigskip

\noindent $\bullet$ We prove (c). Let $t\in A$, so that $ts_{1}^\infty$, 
$ts_{2}s_{1}^\infty\in A^\infty$. These sequences split after $t(s_{1}s_{2}\wedge 
s_{2}s_{1})$, and the decomposition of $ts_{1}^\infty$ (resp., $ts_{2}s_{1}^\infty$) 
into words of $\{s_{1},s_{2}\}$ starts with $us_{i}$ (resp., $us_{3-i}$), 
where ${u\in [\{s_{1},s_{2}\}^{<\omega}]^{*}}$. So $ts_{1}^\infty$ and 
$ts_{2}s_{1}^\infty$ split after $u(s_{1}s_{2}\wedge s_{2}s_{1})$ by (b). 
But we must have $t=u$ because of the position of the splitting point.\bigskip

\noindent $\bullet$ We prove Lemma 8. If $A\in {\cal G}_{0}$, then $F:=\emptyset$ 
works. If $A\in {\cal G}_{1}\setminus {\cal G}_{0}$, then let $w\in n^{<\omega}
\setminus\{\emptyset\}$ such that $A\subseteq\{ w^q/q\in\omega\}$, and $q>0$ such that 
$w^q\in A$. Then $F:=\{ w^q\}$ works. So we may assume that $A\notin 
{\cal G}_{1}$, and $A^\infty =\{ s_{1},s_{2}\}^\infty$. As 
$A^\infty\subseteq\bigcup_{t\in A^{-}}~\{\alpha\in N_{t}/
s_{1}s_{2}\wedge s_{2}s_{1}\prec\alpha -t\}$ is compact, we get a finite subset 
$F$ of $A^-$ such that $A^\infty\subseteq\bigcup_{t\in F}~\{\alpha\in N_{t}/
s_{1}s_{2}\wedge s_{2}s_{1}\prec\alpha -t\}$. We have $F^\infty\subseteq 
A^\infty$. If $\alpha\in A^\infty$, then let $t\in F$ such that $t\prec\alpha$. By (c), we have 
$t\in [\{s_{1},s_{2}\}^{<\omega}]^{*}$. The sequence $t$ is the 
beginning of the decomposition of $\alpha$ into words of $\{s_{1},s_{2}\}$. Thus 
$\alpha -t\in A^\infty$ and we can go on like this. This shows that 
$\alpha\in F^\infty$.$\hfill\square$\bigskip

\noindent\bf Remark.\rm ~The inclusion of $A^\infty =\{s_{1},s_{2}\}^\infty$ 
into $\{t_{1},t_{2}\}^\infty$ does not imply 
$\{s_{1},s_{2}\}\subseteq [\{t_{1},t_{2}\}^{<\omega}]^{*}$, even if 
$A\notin {\cal G}_{1}$. Indeed, take $s_{1}:=01$, $s_{2}:=t_{1}:=0$ and $t_{2}:=10$. 
But we have 
$$\vert t_{1}\vert +\vert t_{2}\vert \leq\vert s_{1}\vert +\vert s_{2}\vert\mbox{,}$$ 
which is the case in general:

\begin{lem} Let $A$, $B\notin {\cal G}_{1}$ satisfying 
$A^\infty =\{s_{1},s_{2}\}^\infty\subseteq B^\infty =\{t_{1},t_{2}\}^\infty$. 
Then there is $j\in 2$ such that $\vert t_{1+i}\vert \leq\vert s_{1+[i+j~\mbox{mod}~2]}\vert $ for each $i\in 2$. In particular, 
$\vert t_{1}\vert  +\vert t_{2}\vert \leq\vert s_{1}\vert  +\vert s_{2}\vert $.\end{lem}

\bf\noindent Proof.\rm ~We may assume that $\vert s_{1}\vert \leq
\vert s_{2}\vert $. Let, for $i=1,~2$, $(w^i_{m})_{m}\subseteq\{t_{1},t_{2}\}$ be 
sequences such that $s_{1}^\infty = w^1_{0}w^1_{1}\ldots$ (resp., 
$s_{2}s_{1}^\infty = w^2_{0}w^2_{1}\ldots$). By the proof of Lemma 8, 
there is a minimal integer $m_{0}$ satisfying $w^1_{m_{0}}\not= w^2_{m_{0}}$. We 
let $u:=w^1_{0}\ldots w^1_{m_{0}-1}$. The sequences $s_{1}^\infty$ and 
$s_{2}s_{1}^\infty$ split after $s_{1}s_{2}\wedge s_{2}s_{1}=u
(t_{1}t_{2}\wedge t_{2}t_{1})$. Similarly, $s_{1}^\infty$ and 
$s_{1}s_{2}s_{1}^\infty$ split after $s_{1}(s_{1}s_{2}\wedge s_{2}s_{1})=v
(t_{1}t_{2}\wedge t_{2}t_{1})$, where $v\in [\{ t_{1},t_{2}\}^{<\omega}]^{*}
\setminus\{\emptyset\}$. So we get $s_{1}u=v$. Similarly, with the sequences 
$s_{2}s_{1}^\infty$ and $s_{2}^2s_{1}^\infty$, we see that 
$s_{2}u\in [\{ t_{1},t_{2}\}^{<\omega}]^{*}\setminus\{\emptyset\}$. So we may 
assume that $u\not= \emptyset$ since $\{ s_{1},s_{2}\}\notin {\cal G}_{1}$. If 
$t_{1}\not\perp t_{2}$, then we may assume that $\emptyset\not= 
t_{1}\prec_{\not=}t_{2}$. So we may assume that we are not in the 
case $t_{2}\!\prec\! t_{1}^\infty$. Indeed, otherwise 
$t_{2}=t_{1}^mt$, where $\emptyset\prec_{\not=}t\prec_{\not=}t_{1}$ (see the proof 
of Lemma 8). Moreover, $t_{1}$ doesn't finish $t_{2}$, otherwise we would have 
$t_{1}=t(t_{1}-t)=(t_{1}-t)t$ and $t$, $t_{1}-t$, $t_{1}$, $t_{2}$ 
would be powers of the same sequence, which contradicts $\{t_{1},t_{2}\}\notin 
{\cal G}_{1}$. As $s_{i}u\in [\{ t_{1},t_{2}\}^{<\omega}]^{*}$, 
this shows that $s_{i}\in [\{ t_{1},t_{2}\}^{<\omega}]^{*}$. So we 
are done since $\{ s_{1},s_{2}\}\notin {\cal G}_{1}$ as before.

\vfill\eject

 Assume for example that $t_{2}=w^1_{m_{0}}$. Let $m'$ be maximal with 
$t_{1}^{m'}\prec t_{2}$. Notice that 
$${ut_{1}^{m'}\prec s_{1}s_{2}\prec s_{1}s_{2}s_{1}^\infty}.$$ 
We have ${ut_{2}\prec s_{1}s_{2}s_{1}^\infty}$, 
otherwise we would obtain ${ut_{1}^{m'+1}\prec s_{1}s_{2}s_{1}^\infty\wedge s_{2}
s_{1}^\infty = s_{1}s_{2}\wedge s_{2}s_{1}\prec s_{1}^\infty}$, which is absurd. 
So we get ${\vert t_{2}\vert \leq\vert s_{1}\vert }$ since ${\vert u\vert \! +\!\vert t_{2}\vert \! +\!
\vert t_{1}t_{2}\wedge t_{2}t_{1}\vert \!\leq\!\vert s_{1}\vert \! +\!
\vert s_{1}s_{2}\wedge s_{2}s_{1}\vert }$. Similarly, $\vert t_{1}\vert \leq\vert s_{2}\vert $ since 
$ut_{1}^{m'+1}\prec s_{2}^2s_{1}^\infty$. The argument is similar if 
$t_{2}=w^2_{m_{0}}$ (we get $\vert t_{i}\vert \leq\vert s_{i}\vert $ in this case for $i=1$, 
$2$).$\hfill\square$ 

\begin{cor} ${\cal G}_{2}$ is a $\check 
D_{\omega}(\boraone)\setminus D_{\omega}(\boraone)$ set. In particular, 
${\cal G}_{2}$ is $\check D_{\omega}(\boraone)$-complete.\end{cor}

\bf\noindent Proof.\rm ~We will apply the Hausdorff derivation 
to ${\cal G}\subseteq 2^{n^{<\omega}}$. This means that we define a 
decreasing sequence $(F_{\xi})_{\xi <\omega_{1}}$ of closed subsets of 
$2^{n^{<\omega}}$ as follows:
$$F_{\xi}:=\overline{\left(\bigcap_{\eta <\xi}~F_{\eta}\right)\cap 
{\cal G}}~\mbox{if}~\xi ~\mbox{is~even},~~\overline{\left(\bigcap_{\eta <\xi}~F_{\eta}\right)
\setminus {\cal G}}~\mbox{if}~\xi ~\mbox{is~odd.}$$
Recall that if $\xi$ is even, then $F_{\xi}=\emptyset$ is equivalent to 
${\cal G}\in D_{\xi}(\boraone)$. Indeed, we set $U_{\xi}:=\check F_{\xi}$. We have $U_{\xi +1}\setminus 
U_{\xi}=F_{\xi}\setminus F_{\xi +1}\subseteq {\cal G}$ if $\xi$ is even and 
$U_{\xi +1}\setminus U_{\xi}\subseteq\check {\cal G}$ if $\xi$ is odd. 
Similarly, ${U_{\xi}\setminus (\bigcup_{\eta 
<\xi}~U_{\eta})\subseteq\check {\cal G}}$ if $\xi$ is limit. If 
$F_{\xi}=\emptyset$, then let $\eta$ be minimal such that $F_{\eta}=\emptyset$. We 
have ${{\cal G}=\bigcup_{\theta\leq\eta ,~\theta ~\mbox{odd}}~U_{\theta}\setminus 
(\bigcup_{\rho <\theta}~U_{\rho})}$. If $\eta$ is odd, then ${\check {\cal G}=
\bigcup_{\theta <\eta ,~\theta ~\mbox{even}}~U_{\theta}\setminus 
(\bigcup_{\rho <\theta}~U_{\rho})\in D_{\eta}(\boraone)}$, thus 
${\cal G}\in \check D_{\eta}(\boraone)\subseteq D_{\xi}(\boraone)$. If $\eta$ is 
even, then ${{\cal G}=\bigcup_{\theta <\eta ,~\theta ~\mbox{odd}}~U_{\theta}\setminus 
(\bigcup_{\rho <\theta}~U_{\rho})\in D_{\eta}(\boraone)}$ and the same conclusion is true. Conversely, if 
${\cal G}\in D_{\xi}(\boraone)$, then let $(V_{\eta})_{\eta <\xi}$ be an increasing sequence of open 
sets with ${\cal G}=\bigcup_{\eta <\xi ,~\eta ~\mbox{odd}}~V_{\eta}\setminus 
(\bigcup_{\theta <\eta}~V_{\theta})$. By induction, we check that 
$F_{\eta}\subseteq\check V_{\eta}$ if $\eta <\xi$. This clearly 
implies that $F_{\xi}=\emptyset$ because $\xi$ is even.\bigskip

\noindent $\bullet$ We will show that if 
$A\notin {\cal G}_{1}$ satisfies $A^\infty = \{s_{1},s_{2}\}^\infty$, then 
$A\notin F_{M}:=F_{M}({\cal G}_{2})$, where $M$ is the smallest odd integer 
greater than or equal to $f(s_{1},s_{2}):=2\Sigma_{l\leq\vert s_{1}\vert +
\vert s_{2}\vert -2}~n^{2(\vert s_{1}\vert +\vert s_{2}\vert -l)}$.\bigskip

 We argue by contradiction: $A$ is the limit of $(A_{q})$, where 
$A_{q}\in F_{M-1}\setminus {\cal G}_{2}$. Lemma 8 gives a finite subset $F$ of 
$A$, and we may assume that $F\subseteq A_{q}$ for each $q$. Thus we have 
$A^\infty\subseteq A_{q}^\infty$, and the inclusion is strict. Thus we can find 
$s^q\in [A_{q}^{<\omega}]^{*}$ such that $N_{s^q}\cap A^\infty =\emptyset$. 
Let $s^q_{0}$, $\ldots$, $s^q_{m_{q}}\in A_{q}$ be such that $s^q=s^q_{0}\ldots 
s^q_{m_{q}}$.\bigskip 

 Now $A_{q}$ is the limit of $(A_{q,r})_{r}$, where 
${A_{q,r}\in F_{M-2}\cap {\cal G}_{2}}$, and we may assume that 
$${\{s^q_{0},~\ldots,~s^q_{m_{q}}\}\cup F\subseteq A_{q,r}}$$ 
for each $r$, and that $A_{q,r}\notin {\cal G}_{1}$ because $A_{q}\notin {\cal G}_{1}\subseteq 
{\cal G}_{2}$. Let $s_{1}^{q,r}$, $s_{2}^{q,r}$ such that $A_{q,r}^\infty =
\{s_{1}^{q,r}, s_{2}^{q,r}\}^\infty$. By Lemma 9 we have $\vert s_{1}^{q,r}\vert +
\vert s_{2}^{q,r}\vert \leq\vert s_{1}\vert +\vert s_{2}\vert $. Now we let $B_{0}:=A_{0,0}$ and 
$s_{i}^0:=s_{i}^{0,0}$ for $i=1,~2$. We have ${B_{0}\in F_{M-2}\cap {\cal G}_{2}
\setminus {\cal G}_{1}}$, $A^\infty\subseteq_{\not=}B_{0}^\infty =
\{s_{1}^0,s_{2}^0\}^\infty$, and 
$$\vert s_{1}^{0}\vert +\vert s_{2}^{0}\vert \leq\vert s_{1}\vert +\vert s_{2}\vert .$$ 

\vfill\eject

 Now we iterate this: for each ${0<k<n^{2(\vert s_{1}\vert +\vert s_{2}\vert )}}$, we get 
${B_{k}\in F_{M-2(k+1)}\cap {\cal G}_{2}\setminus {\cal G}_{1}}$ such that 
${B_{k-1}^\infty\subseteq_{\not=}B_{k}^\infty =\{s_{1}^k,s_{2}^k\}^\infty}$ and 
${\vert s_{1}^{k}\vert +\vert s_{2}^{k}\vert \leq\vert s_{1}^{k-1}\vert +\vert s_{2}^{k-1}\vert }$. We 
can find $k_{0}<n^{2(\vert s_{1}\vert +\vert s_{2}\vert )}$ such that 
${\vert s_{1}^{k_{0}}\vert +\vert s_{2}^{k_{0}}\vert <\vert s_{1}^{k_{0}-1}\vert +
\vert s_{2}^{k_{0}-1}\vert }$ (with the convention $s_{i}^{-1}:=s_{i}$). We set $C_{0}:=
B_{k_{0}}$, $t_{i}^0:=s_{i}^{k_{0}}$. So we have $C_{0}\in F_{M-2(k_{0}+1)}
\cap {\cal G}_{2}\setminus {\cal G}_{1}$, $C_{0}^\infty =\{t_{1}^0,t_{2}^0\}^
\infty$ and ${\vert t_{1}^{0}\vert +\vert t_{2}^{0}\vert <\vert s_{1}\vert +\vert s_{2}\vert }$. Now we 
iterate this: for each $l\leq\vert s_{1}\vert +\vert s_{2}\vert -2$, we get $t_{1}^{l}$, 
$t_{2}^{l}$, $k_{l}<n^{2(\vert t_{1}^{l-1}\vert +\vert t_{2}^{l-1}\vert )}$ and 
$$C_{l}\in F_{M-2\Sigma_{m\leq l}~(k_{m}+1)}\cap {\cal G}_{2}\setminus {\cal G}_{1}$$ 
satisfying $C_{l}^\infty =\{t_{1}^{l},t_{2}^{l}\}^\infty$ and 
${\vert t_{1}^{l}\vert +\vert t_{2}^{l}\vert <\vert t_{1}^{l-1}\vert +\vert t_{2}^{l-1}\vert }$ 
(with the convention $t_{i}^{-1}:=s_{i}$). We have 
${\vert t_{1}^{l}\vert +\vert t_{2}^{l}\vert \leq \vert s_{1}\vert +\vert s_{2}\vert -1-l}$, thus 
$${2\Sigma_{l\leq\vert s_{1}\vert +\vert s_{2}\vert -2}~(k_{l}+1)\leq 
2\Sigma_{l\leq\vert s_{1}\vert +\vert s_{2}\vert -2}~n^{2(\vert t_{1}^{l-1}\vert +
\vert t_{2}^{l-1}\vert )}\leq f(s_{1},s_{2})}$$ 
and this construction is possible. But we have $\vert t_{1}^{\vert s_{1}\vert +\vert s_{2}\vert -2}\vert +
\vert t_{2}^{\vert s_{1}\vert +\vert s_{2}\vert -2}\vert \leq 1$, thus 
${C_{\vert s_{1}\vert +\vert s_{2}\vert -2}\in {\cal G}_{1}}$, which is absurd.\bigskip

 Let $A\notin {\cal G}_{2}$. As $A\notin {\cal G}_{1}$, we can find $s,~t\in A$ 
which are not powers of the same sequence. Indeed, let $s\in A^-$ and $u$ with 
minimal length such that $s$ is a power of $u$. Then any 
$t\in A\setminus\{u^q/q\in\omega\}$ works, because if $s$ and $t$ are powers of 
$w$, then $w$ has to be a power of $u$. Indeed, as $u\prec w$, $w=u^kv$ with 
$v\prec u$, and $v$ has to be a power of $u$ by minimality of $\vert u\vert $ and Lemma 
7. Assume that moreover $A\in F_{2k+2}$. Now $A$ is the limit of 
$(A_{k,r})_{r}\subseteq F_{2k+1}\cap {\cal G}_{2}$ for each integer $k$, and we 
may assume that $s,~t\in A_{k,r}\notin {\cal G}_{1}$. Let $s_{1}^{k,r}$, 
$s_{2}^{k,r}$ be such that ${A_{k,r}}^\infty=\{s_{1}^{k,r},s_{2}^{k,r}\}^\infty$. 
By Lemma 9 we have $\vert s_{1}^{k,r}\vert +\vert s_{2}^{k,r}\vert \leq \vert s\vert +\vert t\vert $ 
and $f(s_{1}^{k,r},s_{2}^{k,r})\leq f(s,t)$. By the preceding point, we must have 
$$2k+1<f(s,t).$$ 
Thus $\bigcap_{m}~F_{m}\subseteq {\cal G}_{2}$. Notice that 
$F_{m+1}(\check {\cal G}_{2})\subseteq F_{m}$, so that $F_{\omega}(\check {\cal G}_{2})=\emptyset$ and 
${\cal G}_{2}\in\check D_{\omega}(\boraone)$.\bigskip

\noindent $\bullet$ Now let us show that $\{0\}\in F_{\omega}({\cal G}_{2})$  
$($this will imply ${\cal G}_{2}\notin D_{\omega}(\boraone))$. It is enough to see that 
$$\{0\}\in \bigcap_{m}~F_{m}.$$
Let $E(x)$ be the biggest integer less than or equal to $x$, 
${p_{k,s}:=2^{k+1-E(\vert s\vert /2)}}$ and 
$k\in\omega$. We define $A_{\emptyset}:=\{0\}$ and, for $s\in 
(\omega\setminus\{0,1\})^{\leq 2k+1}$ and 
$m>1$, ${A_{sm}\! :=\! A_{s}\! \cup\!  \{ {(01^{ p_{k,s} })}^m;
{({0^2}1^{ p_{k,s} })}^m\}}$ if $\vert s\vert $ is even, 
${A_{s}\! \cup\!  \{s\in [\{0,1^{p_{k,s}}\}^{<\omega}]^{*}\! /m\! \leq\! 
\vert s\vert \! \leq\!  m\! +\! 
p_{k,s}\}}$ if $\vert s\vert $ is odd. Let us show that $A_{s}\in {\cal 
G}_{2}$ (resp., $\check {\cal G}_{2}$) if $\vert s\vert $ is even (resp., odd). 
First by induction we get ${A_{sm}\subseteq 
\{0,1^{p_{k,s} }\}^{<\omega}}$. Therefore $A_{sm}^\infty =
\{0,1^{p_{k,s} }\}^\infty$ if $\vert s\vert $ is odd, because 
if $\alpha$ is in $\{0,1^{p_{k,s} }\}^\infty$ and  
$t\in [\{0,1^{p_{k,s} }\}^{<\omega}]^{*}$ with minimal length $\geq m$ begins 
$\alpha$, then ${t\in A_{sm}}$. Now if $\vert s\vert $ is 
even and ${A_{sm}^\infty=\{s_{1},s_{2}\}^\infty}$, then $0^\infty\in\{s_{1},s_{2}
\}^\infty$, thus for example $s_{1}=0^{k+1}$. ${{(01^{ p_{k,s}  
})}^\infty\in \{s_{1},s_{2}\}^\infty}$, thus $s_{2}\prec {(01^{p_{k,s}})}^\infty$ 
and ${\vert s_{2}\vert \geq\vert {(01^{ p_{k,s} })}^m\vert }$ since $s_{2}0^\infty\in
\{s_{1},s_{2}\}^\infty$. But then ${({0^2}1^{ p_{k,s} })}^\infty\notin
\{s_{1},s_{2}\}^\infty$ since $m>1$. 
Thus $A_{sm}\notin {\cal G}_{2}$.\bigskip

 As $(A_{sm})_{m}$ tends to $A_{s}$ and $(A_{s})_{\vert s\vert =2k+2}\subseteq 
{\cal G}_{2}$, we deduce from this that ${A_{s}}$ is in ${F_{2k+1-\vert s\vert }\setminus 
{\cal G}_{2}}$ if $\vert s\vert \leq 2k+1$ is odd, and that $A_{s}\in F_{2k+1-\vert s\vert }\cap 
{\cal G}_{2}$ if $\vert s\vert \leq 2k+1$ is even. Therefore $\{0\}$ is in $\bigcap_{k}~
F_{2k+1}=\bigcap_{m}~F_{m}$.$\hfill\square$

\vfill\eject

\noindent\bf Remarks.\rm ~(1) The end of this proof also shows that ${\cal 
G}_{p}\notin D_{\omega}(\boraone)$ if $p\geq 2$. Indeed, 
$\{ 0\}\in F_{\omega}({\cal G}_{p})$. The only thing to change is the 
definition of $A_{sm}$ if $\vert s\vert $ is even: we set 
$${A_{sm}:=A_{s}\cup\{ (0^{j+1}1^{p_{k,s}})^m/j<p\}}.$$
(2) If $\{s_{1},s_{2}\}\notin {\cal G}_{1}$ and
$\{s_{1},s_{2}\}^\infty =\{t_{1},t_{2}\}^\infty$, then $\{s_{1},s_{2}\}=
\{t_{1},t_{2}\}$. Indeed, $\{t_{1},t_{2}\}\notin {\cal G}_{1}$, thus 
by Lemma 9 we get $\vert s_{1}\vert +\vert s_{2}\vert =\vert t_{1}\vert +\vert t_{2}\vert $. 
By (c) in the proof of Lemma 8 and the previous fact, $s_{i}=t_{\varepsilon_{i}}^{a_{i}}$, where 
$a_{i}>0$, $\varepsilon_{i},~i\in\{ 1,2\}$. As $\{s_{1},s_{2}\}\notin {\cal G}_{1}$, 
$\varepsilon_{1}\not=\varepsilon_{2}$. Thus $a_{i}=1$.\bigskip

\noindent\bf Conjecture 1.\rm ~Let $A\in {\cal F}$. Then there exists a finite 
subset $F$ of $A$ such that $A^\infty=F^\infty$.\bigskip

\noindent\bf Conjecture 2.\rm ~Let $p\geq 1$, $A$, $B\notin {\cal G}_{p}$ 
with ${A^\infty =\{s_{1},\ldots ,s_{q}\}^\infty\subseteq B^\infty =
\{t_{1},\ldots ,t_{p+1}\}^\infty}$. Then 
$\Sigma_{1\leq i\leq p+1}~\vert t_{i}\vert \leq\Sigma_{1\leq i\leq q}~\vert s_{i}\vert $.\bigskip

\noindent\bf Conjecture 3.\rm ~We have ${\cal G}_{p+1}\setminus {\cal 
G}_{p}\in D_{\omega}(\boraone)$ for each $p\geq 1$. In particular, 
${\cal F}\in K_{\sigma}\setminus\bormtwo$.\bigskip

 Notice that Conjectures 1 and 2 imply Conjecture 3. Indeed, 
${\cal F}={\cal G}_{1}\cup\bigcup_{p\geq 1}~{\cal G}_{p+1}\setminus {\cal G}_{p}$, 
so ${\cal F}\in K_{\sigma}$ if ${\cal G}_{p+1}\setminus {\cal G}_{p}\in 
D_{\omega}(\boraone)\subseteq\bortwo$, by Proposition 6. By Proposition 5 
we have ${\cal F}\notin\bormtwo$. It is enough to see that 
$F_{\omega}:=F_{\omega}
({\cal G}_{p+1}\setminus {\cal G}_{p})=\emptyset$. We argue as in the 
proof of Corollary 10. This time, 
$f(s_{1},\ldots ,s_{q}):=2\Sigma_{l\leq\Sigma_{1\leq i\leq q}~\vert s_{i}\vert -2}~n^{q(\Sigma_{1\leq i\leq q}~\vert s_{i}\vert -l)}$ for ${s_{1},~\ldots ,~s_{q}\in n^{<\omega}}$. The 
fact to notice is that $A\notin F_{M}({\cal G}_{p+1}\setminus {\cal G}_{p})$ if 
$A\notin {\cal G}_{p}$ satisfies ${A^\infty =\{s_{1},\ldots ,s_{p+1}\}^\infty}$ and $M$ is 
the minimal odd integer greater than or equal to $f(s_{1},\ldots 
,s_{p+1})$. So if ${A\in F_{2k+2}\cap {\cal F}\setminus {\cal G}_{p}}$, then 
Conjecture 1 gives a finite subset $F:=\{s_{1},\ldots ,s_{q}\}$ of $A$. The set $A$ is 
the limit of 
${(A_{k,r})_{r}\subseteq F_{2k+1}\cap {\cal G}_{p+1}\setminus {\cal G}_{p}}$ for 
each integer $k$, and we may assume that $F\subseteq A_{k,r}$. Conjecture 
2 implies that $f(s_{1}^{k,r},\ldots ,s_{p+1}^{k,r})\leq 
f(s_{1},\ldots ,s_{q})$ and ${2k+1<f(s_{1},\ldots ,s_{q})}$. Thus 
$\bigcap_{m}~F_{m}\subseteq\check {\cal F}\cup {\cal G}_{p}$. So 
$F_{\omega}\subseteq\overline{(\check {\cal F}\cup {\cal G}_{p})\cap 
{\cal G}_{p+1}\setminus {\cal G}_{p}}=\emptyset$.\bigskip\smallskip 

\noindent\bf {\Large 3 Is $A^\infty$ Borel?}\bigskip\rm

 Now we will see that the maximal complexity is possible. We essentially 
give O. Finkel's example, in a lightly simpler version.
 
\begin{prop} Let $\Gamma:=\ana$ or a 
Baire class. The existence of $n\in\omega\setminus 2$ and  
$A\subseteq n^{<\omega}$ such that $A^\infty$ is $\Gamma$-complete is 
equivalent to the existence of $B\subseteq 2^{<\omega}$ such that $B^\infty$ 
is $\Gamma$-complete.\end{prop}

\bf\noindent Proof.\rm ~Let $p_{n}:=\mbox{min}\{p\in\omega/n\leq 
2^p\}\geq 1$. We define $\phi : n\hookrightarrow 2^{p_{n}}:=\{\sigma_{0},\ldots ,\sigma_{2^{p_{n}}-1}\}$ by the formula $\phi (m):=\sigma_{m}$, $\Phi : n^{<\omega}
\hookrightarrow 2^{<\omega}$ by the formula $\Phi (t):=\phi\big(t(0)\big)\ldots 
\phi\big(t(\vert t\vert -1)\big)$ and $f : n^\omega\hookrightarrow 2^\omega$ 
by the formula $f(\gamma ):=\phi\big(\gamma (0)\big)\phi\big(\gamma (1)\big)\ldots$ 
Then $f$ is an homeomorphism from $n^\omega$ onto its range and reduces $A^\infty$ 
to $B^\infty$, where $B := \Phi [A]$. The inverse function of $f$ 
reduces $B^\infty$ to $A^\infty$. So we are done if $\Gamma$ is 
stable under intersection with closed sets. Otherwise, $\Gamma=\borone$ or 
$\boraone$. If $A=\{s\in 2^{<\omega}/0\prec s~\mbox{or}~1^2\prec s\}$, then 
$A^\infty =N_{0}\cup N_{1^2}$, which is $\borone$-complete. If 
${A = \{s\in 2^{<\omega}/0\prec s\}\cup\{10^k1^{l+1}/k,l\in
\omega\}}$, then ${A^\infty =2^\omega\setminus\{10^\infty\}}$, which is 
$\boraone$-complete.$\hfill\square$

\begin{thm} The set $I:=\{(\alpha ,A)\in n^\omega
\times 2^{n^{<\omega}}/\alpha\in A^\infty\}$ is $\ana$-complete. In fact,\smallskip

\noindent (a) (O. Finkel, see~[F1]) There exists $A_{0}\subseteq 2^{<\omega}$ such that $A_{0}^\infty$ is 
$\ana$-complete.\smallskip 

\noindent (b) There exists $\alpha_{0}\in2^\omega$ such that $I_{\alpha_{0}}$ is 
$\ana$-complete.\end{thm}

\vfill\eject

\bf\noindent Proof.\rm ~(a) We set $L:=\{2,3\}$ and\bigskip
 
\leftline{${\cal T} := \{~\tau\subseteq 2^{<\omega}\times L/ \forall (u,\nu )\in 2^{<\omega}\times L~~
[(u,\nu )\notin\tau ]~\mbox{or}$}\bigskip

\rightline{${[\big(\forall v\prec u~\exists \mu\in L~(v,\mu )\in\tau \big)~\mbox{and}~\big((u,5-\nu )\notin\tau\big)~\mbox{and}~\big(\exists 
 (\varepsilon,\pi )\in 2\times L~(u\varepsilon ,\pi )\in\tau\big)]~\}}$.}\bigskip
 
\noindent The set $\cal T$ is the set of pruned trees over $2$ with labels in $L$. It is 
a closed subset of $2^{2^{<\omega}\times L}$, thus a Polish space. Then we set 
$$\sigma := \{\tau\in {\cal T}/\exists (\underline{u},\underline{\nu})\in 2^\omega\times L^\omega ~
[\forall m~\big(\underline{u}\lceil m,\underline{\nu }(m)\big)\in\tau ]~\mbox{and}~[\forall p~
\exists m\geq p~~\underline{\nu }(m)=3]\}.$$
$\bullet$ Then $\sigma\in\ana ({\cal T})$. Let us show that it is 
complete. We set ${{\cal T}:=\{T\in 2^{\omega ^{<\omega}}/T~\mbox{is~a~tree}\}}$ and 
${IF := \{T\in {\cal T}/T~\mbox{is~ill-founded}\}}$. 
It is a well-known fact that ${\cal T}$ is a Polish space (it is a 
closed subset of $2^{\omega ^{<\omega}}$), and that $IF$ is $\ana$-complete 
(see [K1]). It is enough to find a Borel reduction of $IF$ to $\sigma$ (see [K2]).\bigskip 

We define $\psi :\omega^{<\omega}\hookrightarrow 2^{<\omega}$ by the 
formula $\psi (t):= 0^{t(0)}10^{t(1)}1\ldots 0^{t(\vert t\vert -1)}1$, and 
$\Psi :{\cal T}\rightarrow {\cal T}$ by\bigskip

\leftline{$\Psi (T):=\{ (u,\nu )\in 2^{<\omega}\times L/\exists t\in T~u\prec\psi (t)~\mbox{and}~\nu =  3~
\mbox{if}~u=\emptyset ,~2+u(\vert u\vert -1)~\mbox{otherwise} \}$}\bigskip

\rightline{$\cup ~\{ (\psi (t)0^{k+1},2)/t\in T~\mbox{and}~\forall q\in\omega~t 
q\notin T, k\in\omega \}.$}\bigskip

\noindent The map $\Psi$ is Baire class one. Let us show that it is a reduction. If $T\in IF$, then 
let $\gamma\in\omega^\omega$ be such that $\gamma\lceil m\in T$ 
for each integer $m$. We have $\big(\psi (\gamma\lceil m ),3\big)\in\Psi (T)$. 
Let $\underline{u}$ be the limit of $\psi (\gamma\lceil m)$ and  
$\underline{\nu }(m):= 2+\underline{u} (m-1)$ (resp., $3$) if $m>0$ 
(resp., $m=0$). These objects show that $\Psi (T)\in\sigma$. Conversely, we have 
$T\in IF$ if $\Psi (T)\in\sigma$.\bigskip

\noindent $\bullet$ If $\tau\in {\cal T}$ and $m\in\omega$, then we enumerate 
$\tau\cap (2^m\times L) := \{(u^{m,\tau }_{1},\nu ^{m,\tau }_{1}),\ldots ,(u^{m,
\tau }_{q_{m,\tau }},\nu ^{m,\tau }_{q_{m,\tau }})\}$ in the lexicographic 
ordering. We define $\varphi : {\cal T}\hookrightarrow 5^\omega$ by the formula 
$$\varphi (\tau ):=
( {u^{0,\tau }_{1}} {\nu ^{0,\tau }_{1}}\ldots 
  {u^{0,\tau }_{q_{0,\tau }}} {\nu ^{0,\tau }_{q_{0,\tau }}}4) 
( {u^{1,\tau }_{1}} {\nu ^{1,\tau }_{1}}\ldots 
  {u^{1,\tau }_{q_{1,\tau }}} {\nu ^{1,\tau }_{q_{1,\tau 
  }}}4)\ldots$$
The set $A_{0}$ will be made of finite subsequences of sentences in $\varphi [{\cal T}]$. We set 
$$\begin{array}{ll}
A_{0}:=\{\!\!\!\! 
& {{u^{m,\tau }_{q+1}}}{\nu ^{m,\tau }_{q+1}}\ldots {
{u^{p,\tau }_{r}}}{\nu ^{p,\tau }_{r}}/\tau\in 
{\cal T},~m+1\! <\! p,~0\leq q\leq q_{m,\tau},~1\! \leq\!  r\! \leq\!  q_{p,\tau},\cr 
& [(m=0~\mbox{and}~q=0)~\mbox{or}~(q>0~\mbox{and}~{\nu ^{m,\tau 
 }_{q}}={3}~\mbox{and}~{u ^{m,\tau }_{q}}\prec 
{u ^{p,\tau }_{r}})],~{\nu ^{p,\tau }_{r}}=3\ \}
\end{array}$$
(with the convention ${{u^{m,\tau }_{q_{m,\tau }+1}}}{\nu ^{m,\tau 
}_{q_{m,\tau }+1}}=4$). It is clear that $\varphi$ is continuous, and it 
is enough to see that it reduces $\sigma$ to $A_{0}^\infty$.\bigskip

  So let us assume that $\tau\in\sigma$. This means the existence of 
an infinite branch in the tree with infinitely many $3$ labels. We 
cut $\varphi (\tau )$ after the first $3$ label of the branch 
corresponding to a sequence of length $m>1$. Then we cut after the first 
$3$ label corresponding to a sequence of length at least $m+2$ of the 
branch. And so on. This clearly gives a decomposition of $\varphi (\tau )$ 
into words in $A_{0}$.\bigskip

 If such a decomposition exists, then the first word is ${{u^{0,\tau }_{1}}}
 {\nu ^{0,\tau }_{1}}\ldots {
{u^{p_{0},\tau }_{r_{0}}}}{\nu ^{p_{0},
\tau }_{r_{0}}}$, and the second is ${{{u^{p_{0},\tau }_{r_{0}+1}}}
{\nu ^{p_{0},\tau }_{r_{0}+1}}\ldots {
{u^{p_{1},\tau }_{r_{1}}}}{\nu ^{p_{1},\tau 
}_{r_{1}}}}$. So we have ${u ^{p_{0},\tau }_{r_{0}}}\prec_{\not=} 
{u ^{p_{1},\tau }_{r_{1}}}$. And so on. This gives an infinite 
branch with infinitely many $3$ labels.

\vfill\eject

\noindent $\bullet$ By Proposition 11, we can also have $A_{0}\subseteq 
2^{<\omega}$.\bigskip 

\noindent (b) Let $\alpha_{0}:=1010^210^3\ldots$, $(q_{l})$ be the 
sequence of prime numbers: $q_{0}:=2$, $q_{1}:=3$, $M :\omega^{<\omega}\rightarrow\omega$ 
defined by $M_{s} := q_{0}^{s(0)+1}\ldots q_{\vert s\vert -1}^{s(\vert s\vert -1)+1} +1$, 
$\phi :\omega^{<\omega}\rightarrow 2^{<\omega}\setminus\{\emptyset\}$ defined by the formulas 
$$\phi (\emptyset):=1010^2=1010^{2M_{\emptyset}}$$ 
and $\phi (sm):=10^{2M_{s}+1}10^{2M_{s}+2}\ldots 10^{2M_{sm}}$, and $\Phi : 2^{
\omega ^{<\omega}}\rightarrow 2^{n^{<\omega}}$ defined by $\Phi (T):=\phi [T]$.\bigskip 

\noindent $\bullet$ It is clear that $M_{sm}>M_{s}$, and that $M$ and $\phi$ 
are well defined and one-to-one. So $\Phi$ is continuous:
$$\begin{array}{ll}
s\in\phi [T]\!\!\!\! 
& \Leftrightarrow \exists t~~(t\in T~\mbox{and}~\phi (t)=s)\cr 
& \Leftrightarrow s\in\phi[\omega^{<\omega}]~\mbox{and}~\forall 
t~~(t\in T~\mbox{or}~\phi (t)\not= s).
\end{array}$$
If $T\in IF$, then we can find $\beta\in\omega^\omega$ such that 
$\phi (\beta\lceil l)\in\Phi (T)$ for each integer $l$. Thus 
$$\alpha_{0}=(1010^{2M_{\beta\lceil 0}})(10^{2M_{\beta\lceil 0}+1}\ldots 
10^{2M_{\beta\lceil 1}})\ldots\in\big(\Phi(T)\big)^\infty .$$ 
Conversely, if  $\alpha_{0}\in\big(\Phi(T)\big)^\infty$, then there exist $t_{i}\in T$ such 
that $\alpha_{0}=\phi (t_{0})\phi (t_{1})\ldots$ We have  
$t_{0}=\emptyset$, and, if $i>0$, then $M_{t_{i}\lceil\vert t_{i}\vert -1}=M_{t_{i-1}}$; 
from this we deduce that $t_{i}\lceil\vert t_{i}\vert -1 = t_{i-1}$, because $M$ is 
one-to-one. 
So let $\beta$ be the limit of the $t_{i}$'s. We have $\beta\lceil i=t_{i}$, 
thus $\beta\in [T]$ and $T\in IF$. Thus $\Phi_{\lceil {\cal T}}$ 
reduces $IF$ to $I_{\alpha_{0}}$. Therefore this last set is $\ana$-complete. 
Indeed, it is clear that $I$ is $\Ana$:
$$\alpha\!\in\! A^\infty\Leftrightarrow~\exists\beta\!\in\!\omega^\omega ~~
[\big(\forall m>0~~\beta (m)>0\big)~\mbox{and}~\big(\forall q\!\in\!\omega~~
\pi (\alpha ,\beta ,q)\!\in\! A\big)].$$
Finally, the map from ${\cal T}$ into $n^\omega\times 2^{n^{<\omega}}$, which 
associates $\big(\alpha_{0},\Phi (T)\big)$ to $T$ clearly reduces $IF$ to $I$. So 
$I$ is $\ana$-complete.$\hfill\square$\bigskip

\noindent\bf Remark.\rm ~This proof shows that if $\alpha = 
{s_{0}}{s_{1}}\ldots$ and $(s_{i})$ is an antichain 
for the extension ordering, then $I_{\alpha}$ is $\ana$-complete (here we have 
$s_{i}=10^{2i+1}10^{2i+2})$. To see it, it is enough to notice that 
$\phi (\emptyset)=s_{0}$ and $\phi (sm)={s_{M_{s}}}\ldots 
s_{M_{sm}-1}$. So $I_{\alpha}$ is $\ana$-complete for 
a dense set of $\alpha$'s.\bigskip  

 We will deduce from this some true co-analytic sets. But we need a 
lemma, which has its own interest.

\begin{lem} (a) The set $A^\infty$ is Borel if 
and only if there exist a Borel function $f:n^\omega\rightarrow \omega^\omega$ such that
$$\alpha\in A^\infty ~\Leftrightarrow~\big(\forall m>0~~f(\alpha )(m)>0\big)~\mbox{and}~\big( 
\forall q\in\omega~~\pi (\alpha , f(\alpha ) ,q)\in A\big).$$
(b) Let $\gamma\in\omega^\omega$ and $A\subseteq n^{<\omega}$. Then $A^\infty\in\Borel (A,\gamma )$ if and only if, for $\alpha\in n^\omega$, we have 
$$\alpha\! \in\!  A^\infty\Leftrightarrow~\exists\beta\! \in\! 
\Borel (A,\gamma ,\alpha )~~[\big(\forall m>0~~\beta (m)>0\big)~\mbox{and}~
\big(\forall q\in\omega~~\pi (\alpha ,\beta ,q)\in A\big)].$$\end{lem} 
 
\bf\noindent Proof.\rm ~The ``if" directions in (a) and (b) are clear. We have seen in the proof of 
Proposition 4 the ``if" direction of the equivalences (the existence of an arbitrary 
$\beta$ is necessary and sufficient). So let us show the ``only if" directions.

\vfill\eject

\noindent (a) We define $f:n^\omega\rightarrow\omega^\omega$ by the formula 
$f(\alpha ):=0^\infty ~{\rm if}~\alpha\notin A^\infty$, and, otherwise,
$$f(\alpha )(0):=\mbox{min}\{p\in\omega /\alpha\lceil (p+1)\in A~\mbox{and}~
\alpha -\alpha\lceil (p+1)\in A^\infty\}\mbox{,}$$
\leftline{$f(\alpha )(r+1) :=$}\bigskip

\rightline{$\mbox{min}\{k>0/[\alpha -\alpha\lceil \big(1+\Sigma_{j\leq r}~
f(\alpha )(j)\big)]\lceil k\in A~\mbox{and}~\alpha -\alpha\lceil 
\big( k+1+\Sigma_{j\leq r}~f(\alpha )(j)\big)\in A^\infty\}$.}\bigskip

\noindent We get $\pi (\alpha , f(\alpha ) ,0)= \alpha \lceil f(\alpha )(0)+1\in A$ and, if $q>0$, 
$$\pi (\alpha , f(\alpha ) ,q)= 
~\big(\alpha (1+\Sigma_{j<q}~f(\alpha ) [j]),...,\alpha 
(\Sigma_{j\leq q}~f(\alpha ) [j])\big)~\in A.$$
As $f$ is clearly Borel, we are done.\bigskip

\noindent (b) If $A^\infty\in\Borel (A,\gamma )$, then so is $f$ and $\beta := 
f(\alpha )\in \Borel (A,\gamma ,\alpha )$ is what we were looking 
for.$\hfill\square$\bigskip

\bf\noindent Remark.\rm ~Lemma 13 is a particular case of a more general situation. Actually 
we have the following uniformization result. It was written after a 
conversation with G. Debs.
 
\begin{prop} Let $X$ and $Y$ be Polish spaces, and $F\in\bormtwo (X\times Y)$ such that the projection 
${\it\Pi}_{X}[F\cap (X\times V)]$ is Borel for each $V\in\boraone (Y)$. Then 
there exists a Borel map $f:X\rightarrow Y$ such that 
$\big(x,f(x)\big)\in F$ for each $x\in{\it\Pi}_{X}[F]$.\end{prop}

\bf\noindent Proof.\rm ~Let $(Y_{n})$ be a basis for the 
topology of $Y$ with $Y_{0}:=Y$, $B_{n}:={\it\Pi}_{X}[F\cap (X\times 
Y_{n})]$, and $\tau$ be a finer $0$-dimensional Polish topology on 
$X$ making the $B_{n}$'s clopen (see 13.5 in [K1]). We equip $X$ with a 
complete $\tau$-compatible metric $d$. Let $(O_{m})\subseteq\boraone 
(X\times Y)$ be decreasing satisfying $O_{0}:=X\times Y$ and 
$F=\bigcap_{m}~O_{m}$. We construct a sequence $(U_{s})_{s\in\omega^{<\omega}}$ 
of clopen subsets of $[B_{0},\tau ]$ with $U_{\emptyset}:=B_{0}$, and a sequence 
$(V_{s})_{s\in\omega^{<\omega}}$ of basic open sets of $Y$ satisfying 
$$\begin{array}{ll}
& (a)~ U_{s}\subseteq {\it\Pi}_{X}[F\cap (U_{s}\times V_{s})]\cr
& (b)~\mbox{diam}_{d}(U_{s}),~\mbox{diam}(V_{s})\leq\frac{1}{\vert s\vert }~{\rm 
if}~s\not=\emptyset\cr
& (c)~ U_{s}=\bigcup_{m,\mbox{disj.}}~U_{s^\frown 
m},~\overline{V_{s^\frown n}}\subseteq V_{s}\cr
& (d)~ U_{s}\times V_{s}\subseteq O_{\vert s\vert }
\end{array}$$
$\bullet$ Assume that this construction has been achieved. If 
$x\notin B_{0}$, then we set $f(x):=y_{0}\in Y$ (we may assume that 
$F\not=\emptyset$). Otherwise, we can find a unique sequence $\gamma\in\omega^\omega$ 
such that $x\in U_{\gamma\lceil m}$ for each integer $m$. Thus we can find $y\in 
V_{\gamma\lceil m}$ such that $(x,y)\in F$, and 
$(\overline{V_{\gamma\lceil m}})_{m}$ is a decreasing sequence of 
nonempty closed sets whose diameters tend to $0$, which defines a 
continuous map $f:[B_{0},\tau ]\rightarrow Y$. If $x\in B_{0}$, then 
$\big(x,f(x)\big)\in U_{\gamma\lceil m}\times V_{\gamma\lceil 
m}\subseteq O_{m}$, thus $Gr(f_{\vert B_{0}})\subseteq F$. Notice that 
$f:[X,\tau ]\rightarrow Y$ is continuous, so $f:X\rightarrow Y$ 
is Borel.\bigskip

\noindent $\bullet$ Let us show that the construction is possible. We set 
$U_{\emptyset}:=B_{0}$ and $V_{\emptyset}:=Y$. Assume that 
$(U_{s})_{s\in\omega^{\leq p}}$ and $(V_{s})_{s\in\omega^{\leq p}}$ 
satisfying conditions (a)-(d) have been constructed, which is the 
case for $p=0$. Let $s\in \omega^p$. If $(x,y)\in F\cap (U_{s}\times 
V_{s})$, then we can find $U_{x}\in\borone(U_{s})$ and a basic open set 
$V_{y}\subseteq Y$ such that ${(x,y)\in U_{x}\times V_{y}\subseteq 
U_{x}\times\overline{V_{y}}\subseteq (U_{s}\times V_{s})\cap 
O_{p+1}}$, and whose diameters are at most $\frac{1}{p+1}$. By the Lindel\" of 
property, we can write $F\cap (U_{s}\times V_{s})\subseteq 
\bigcup_{n}~U_{x_{n}}\times V_{y_{n}}$ and $F\cap (U_{s}\times V_{s})=
\bigcup_{n}~F\cap (U_{x_{n}}\times V_{y_{n}})$.

\vfill\eject

 If $x\in U_{s}$, then let 
$n$ and $y$ be such that $(x,y)\in F\cap (U_{x_{n}}\times V_{y_{n}})$. 
Then 
$${x\in O^n:={\it\Pi}_{X}[F\cap (X\times V_{y_{n}})]\cap 
U_{x_{n}}\in\borone([B_{0},\tau ])}.$$ 
Thus $U_{s}=\bigcup_{n}~O^n$. We set 
$U_{s^\frown n}:=O^n\setminus (\bigcup_{p<n}~O^p)$ and $V_{s^\frown n}:=
V_{y_{n}}$, and we are done.$\hfill\square$\bigskip

 In our context, ${F=\{(\alpha ,\beta )\! \in\!  n^\omega\! \times\! \omega^\omega /
\big(\forall m\! >\! 0~~\beta (m)\! >\! 0\big)~\mbox{and}~
\big(\forall q\! \in\! \omega~~\pi (\alpha ,\beta ,q)\! \in\!  A\big)\}}$, which is a 
closed subset of $X\times Y$. The projection $\Pi_{X}[F\cap (X\times 
N_{s})]$ is Borel if $A^\infty$ is Borel, since it is 
$\{{S^{*}}\gamma/S\in (A\cap 
n^{s(0)+1})\times\Pi_{0<j<\vert s\vert }~(A\cap 
n^{s(j)})~\mbox{and}~\gamma\in A^\infty\}$.

\begin{thm} The following sets are $\Ca\setminus\borel$:\smallskip

\noindent (a) ${\it\Pi}:=\{(A,\gamma ,\theta )\in 2^{n^{<\omega}}\times 
 \omega^\omega\times\omega^\omega /\theta\in\mbox{WO~and}~A^\infty\in
 {\bf\Pi}^0_{\vert \theta\vert }\cap\Borel (A,\gamma )\}$. The same thing is true with   
 ${\it\Sigma}:=\{(A,\gamma ,\theta )\in 2^{n^{<\omega}}\times 
 \omega^\omega\times\omega^\omega /\theta\in\mbox{WO~and}~A^\infty\in
 {\bf\Sigma}^0_{\vert \theta\vert }\cap\Borel (A,\gamma )\}$.\smallskip
 
\noindent (b) ${\it{\it\Sigma}}_{1}:=\{A\in 2^{n^{<\omega}}/A^\infty\in\boraone
 \cap\Borel (A)\}$. In fact, ${{\it{\it\Sigma}}_{\xi}:=
 \{A\in 2^{n^{<\omega}}/A^\infty\in\boraxi 
 \cap\Borel (A)\}}$ is $\ca\setminus\borel$ if $1\leq\xi<\omega_{1}$. 
 Similarly, ${{\it{\it\Pi}}_{\xi}:=
 \{A\in 2^{n^{<\omega}}/A^\infty\in\bormxi 
 \cap\Borel (A)\}}$ is $\ca\setminus\borel$ if $2\leq\xi<\omega_{1}$.\smallskip
 
\noindent (c) ${\it{\it\Delta}} :=\{A\in 2^{n^{<\omega}}/A^\infty\in\Borel (A)\}$.\end{thm}
 
\bf\noindent Proof.\rm ~Consider the way of coding the Borel sets 
used in [Lou]. By Lemma 13 we get 
$$(A,\gamma ,\theta )\!\in\! {\it{\it\Pi}}\Leftrightarrow \left\{\!\!\!\!\!\!\!\!
\begin{array}{ll} 
& ~\exists p\!\in\!\omega~~P(p,A,\gamma ,\theta )~\mbox{and}~\forall\alpha\!\in\! n^\omega\cr 
& ~\big(\alpha\!\notin\! A^\infty ~\mbox{or}~(p,A,\gamma ,\alpha)\!\in\! C
\big)~\mbox{and}~\big([(p,A,\gamma )\!\in\! W~\mbox{and}~(p,A,\gamma 
,\alpha)\!\notin\! C]~\mbox{or}\cr 
& ~\exists\beta\!\in\!\Borel (A,\gamma ,
\alpha )~~ [(\forall m\! >\! 0~~\beta (m)\! >\! 0)~\mbox{and}~(\forall q\!\in
\!\omega~~\pi (\alpha ,\beta ,q)\!\in\! A)]\big).
\end{array}
\right.$$
This shows that ${\it{\it\Pi}}$ is $\Ca$. The same 
argument works with ${\it{\it\Sigma}}$. From this we can deduce that ${\it\Sigma}_{1}$ is $\Ca$, 
if we forget $\gamma$ and take the section of ${\it\Sigma}$ at 
$\theta\in {\rm WO}\cap\Borel$ such that $\vert \theta\vert =1$. Similarly, 
${\it\Sigma}_{\xi}$ and ${\it\Pi}_{\xi}$ are co-analytic if $\xi\!\geq\! 1$. 
Forgetting $\theta$, we see that the relation ``$A^\infty\!\in\!\Borel (A,\gamma )$" is 
$\Ca$.\bigskip 

\noindent $\bullet$ Let us look at the proof of Theorem 12. We will 
show that if $\xi\geq 1$ (resp., $\xi\geq 2$), then 
${\it\Sigma}_{\xi}\setminus I_{\alpha_{0}}$ (resp., 
${\it\Pi}_{\xi}\setminus I_{\alpha_{0}}$) is a true co-analytic set. To do 
this, we will reduce $WF$ to ${\it\Sigma}_{\xi}\setminus I_{\alpha_{0}}$ (resp., 
${\it\Pi}_{\xi}\setminus I_{\alpha_{0}}$) in a Borel way. 
We change the definition of $\Phi$. We set 
$$t\subseteq\alpha_{0}~\Leftrightarrow ~\exists 
k~~t\prec\alpha_{0}\! -\!\alpha_{0}\lceil k\mbox{,}$$
$$E := \{ (\alpha_{0}\lceil p)  
r/p\!\in\!\omega\!\setminus\!\{ 2\},r\!\in\! n\!\setminus\!\{\alpha_{0}(p)\}\}\mbox{,}~~~~~~~~~~
F:=\{ U^{*}\!\not\subseteq\!\alpha_{0}/U\!\in\!\phi[T]^{<\omega}\}\mbox{,}$$
$$\Phi'(T):=\phi[T]\cup\{s\!\in\! n^{<\omega}/\exists t\!\in\! E
\!\cup\! F~~t\prec s\}.$$
This time, $\Phi'$ is Baire class one, since 
$$\begin{array}{ll} 
s\!\in\!\Phi'(T)~\Leftrightarrow \!\!
& s\!\in\!\phi [T]~~\mbox{or}~~\exists t\!\in\! E~~t\prec s~~\mbox{or}\cr 
& \exists U\!\in\! (2^{<\omega})^
{<\omega}~~(\forall j\! <\!\vert u\vert ~U(j)\!\in\!\phi[T])~\mbox{and}~U^{*}\!\not
\subseteq\!\alpha_{0}~\mbox{and}~U^{*}\!\prec\! s.
\end{array}$$
The proof of Theorem 12 remains valid, since if $\alpha_{0}\in \big(\Phi'(T)\big)^
\infty$, then the decompositions of $\alpha_{0}$ into words of $\Phi'(T)$ are 
actually decompositions into words of $\phi [T]$.

\vfill\eject

\noindent $\bullet$ Let us show that $\big(\Phi'(T)\big)^\infty\!\in\!\boraone
\cap\Borel \big(\Phi'(T)\big)$ if $T\!\in\! WF$. The set $\big(\Phi'(T)\big)^\infty$ is 
$$\bigcup_{S\in\phi 
[T]^{<\omega},l\in n\setminus\{ 1\},m\in n\setminus\{ 0\}}~
[(\bigcup_{s/\exists t\in F~t\prec s}~N_{S^{*}s})\cup 
N_{S^{*}l}\cup N_{S^{*}1m}\cup (N_{S^{*}101}\setminus\{ 
S^{*}\alpha_{0}\})].$$
If $\alpha\in n^\omega$, then $\alpha$ contains infinitely many $l\in 
n\setminus\{ 1\}$ or finishes with $1^\infty$. As $1^2$ and the sequences 
beginning with $l$ are in $\Phi'(T)$, the clopen sets are subsets of 
$\big(\Phi'(T)\big)^\infty$ since $\phi [T]$ and the sequences beginning 
with $t\in F$, $l$ or $1m$ are in $\Phi'(T)$. If $\alpha\in N_{S^{*}101}\setminus
\{S^{*}\alpha_{0}\}$, then let $p\geq 3$ be maximal such that $\alpha\lceil 
(\vert S^{*}\vert + p) = {S^{*}} (\alpha_0\lceil p)$. We have  
$\alpha\in \big(\Phi'(T)\big)^\infty$ since the sequences beginning with 
$(\alpha_0\lceil p)r$ are in $\Phi'(T)$. Thus we get  the 
inclusion into $\big(\Phi'(T)\big)^\infty$.\bigskip

 If $\alpha\in \big(\Phi'(T)\big)^\infty$, then $\alpha = {a_{0}}{a_{1}}\ldots$, where $a_{i}\in \Phi'(T)$. Either for all $i$ we have $a_{i}\in\phi [T]$. In this case, there is $i$ such that ${a_{0}} 
\ldots {a_{i}}\not\subseteq\alpha_{0}$, otherwise we could 
find $k$ with $\alpha_{0} -\alpha_{0}\lceil 
k\in\big(\Phi(T)\big)^\infty$. But this contradicts the fact that $T\in 
WF$, as in the proof of Theorem 12. So we have $\alpha\in\bigcup_{\exists t\in 
F~t\prec s} N_{s}$. Or there exists $i$ minimal such that $a_{i}\notin\phi 
[T]$. In this case,\bigskip

\noindent - Either $\exists t\in E~t\prec a_{i}$ and $\alpha\in
\bigcup_{S\in\phi [T]^{<\omega},l\in n\setminus\{ 1\},m\in n\setminus\{ 0\}}~
[N_{S^{*}l}\cup N_{S^{*}1m}\cup (N_{S^{*}101}\setminus\{ S^{*}\alpha_{0}\})]$,\bigskip

\noindent - Or $\exists t\in F~t\prec a_{i}$ and 
$\alpha\in\bigcup_{S\in\phi [T]^{<\omega}}\bigcup_{s/\exists t\in F~t
\prec s}~N_{S^{*}s}$.\bigskip

  From this we deduce that $\big(\Phi'(T)\big)^\infty$ is $\boraone$.\bigskip
  
\noindent  Finally, we have 
$$\alpha\!\in\!\big(\Phi'(T)\big)^\infty \Leftrightarrow~\left\{\!\!\!\!\!\!
\begin{array}{ll} 
& \exists t\!\in\! n^{<\omega}~\exists b\!\in\!\omega^{<\omega}~[
\big(\vert t\vert \! =\! 1\! +\!\Sigma_{j<\vert b\vert }~b(j)\big)~\mbox{and}~ 
\big(\forall 0\! <\! m\! <\! \vert b\vert ~b(m)>0\big)\cr 
& \mbox{and}~\big(\forall q\! <\! \vert b\vert ~\pi (t 0^\infty,b 0^\infty,q)\! \in\!  
\Phi'(T)\big)]~\mbox{and}~[\exists l\!\in\! n\!\setminus\!\{ 1\}~t l\!\prec\! 
\alpha~\mbox{or}~t1^2\!\prec\!\alpha ].
\end{array}\right.$$
This shows that $\big(\Phi'(T)\big)^\infty$ is $\Borel \big(\Phi'(T)\big)$.\bigskip 

 Therefore, $\Phi'_{\lceil {\cal T}}$ reduces $WF$ to  
${\it\Sigma}_{\xi}\setminus I_{\alpha_{0}}$ if $\xi\geq 1$, and to 
${\it\Pi}_{\xi}\setminus I_{\alpha_{0}}$ if $\xi\geq 2$. So these sets are 
true co-analytic sets. But ${\it\Sigma}_{1}\cap I_{\alpha_{0}}$ is $\Ca$, 
by Lemma 13. As ${\it\Sigma}_{1}\setminus I_{\alpha_{0}} = 
{\it\Sigma}_{1}\setminus ({\it\Sigma}_{1}\cap I_{\alpha_{0}})$, ${\it\Sigma}_{1}$ is not Borel. 
Thus ${\it\Sigma}$ is not Borel, 
as before. The argument is similar for ${\it\Sigma}_{\xi}$, ${\it\Pi}_{\xi}$ ($\xi\geq 2$) 
and ${\it\Pi}$. And for ${\it\Delta}$ too.$\hfill\square$\bigskip

\noindent\bf Question.\rm ~Does $A^\infty\in\borel$ imply 
$A^\infty\in\Borel (A)$? Probably not. If the answer is positive, 
${\bf\Delta}$, and more generally 
${\bf\Sigma}_{\xi}$ (for $\xi\geq 1$) and 
${\bf\Pi}_{\xi}$ (for $\xi\geq 2$) are 
true co-analytic sets.\bigskip

\noindent\bf Remark.\rm ~In any case, $\bf\Delta$ is ${\it\Sigma}^1_2$ because 
``$A^\infty\! \in\! \borel$" 
is equivalent to ``${\exists\gamma\! \in\! \omega^\omega ~A^\infty\! \in\! \Borel 
(A,\gamma )}$". This argument shows that ${\bf\Sigma}_{\xi}$ and 
${\bf\Pi}_{\xi}$ are ${\it\Sigma}^1_2 (\theta )$, where $\theta\in WO$ satisfies $\vert \theta\vert =\xi$. We can say more about ${\bf\Pi}_{1}$: it is ${\it\Delta}^1_{2}$. 
Indeed, in [St2] we have the following characterization:
$$A^\infty\!\in\!\bormone~\Leftrightarrow ~\forall\alpha\!\in\!  
n^\omega ~[\forall s\!\in\! n^{<\omega}~(s\!\prec\!\alpha\Rightarrow\exists 
S\!\in\! A^{<\omega}~s\!\prec\! S^{*})]\Rightarrow \alpha\!\in\!  
A^\infty .$$
This gives a ${\it\Pi}^1_{2}$ definition of ${\bf\Pi}_{1}$. The same 
fact is true for ${\bf\Sigma}_{1}$: 

\begin{prop} ${\bf\Sigma}_{1}$ and ${\bf\Pi}_{1}$ are co-nowhere dense 
${\it\Delta}^1_{2}\setminus D_{2}(\boraone)$ subsets of $2^{n^{<\omega}}$. If $\xi\geq 2$, then 
${\bf\Sigma}_{\xi}$ and ${\bf\Pi}_{\xi}$ are co-nowhere dense 
${\bf\Sigma}^1_{2}\setminus D_{2}(\boraone)$ subsets of $2^{n^{<\omega}}$. 
${\bf\Delta}$ is a co-nowhere dense ${\it\Sigma}^1_2\setminus D_{2}(\boraone)$ subset of 
$2^{n^{<\omega}}$.\end{prop} 

\vfill\eject

\bf\noindent Proof.\rm ~We have seen that ${\bf\Sigma}_{1}$ is 
${\it\Sigma}^1_2$; it is also ${\it\Pi}^1_{2}$ because 
$$A^\infty\!\in\!\boraone ~\Leftrightarrow ~\forall\alpha\!\in\! 
n^\omega ~~\alpha\!\notin\! A^\infty ~\mbox{or}~\exists s\!\in\! 
n^{<\omega}~[s\!\prec\!\alpha ~\mbox{and}~\forall\beta\!\in\! n^\omega 
~(s\!\not\prec\!\beta ~\mbox{or}~\beta\!\in\! A^\infty )].$$ 
By Proposition 4, ${\bf\Pi}_{0}$ is co-nowhere dense, and it is a 
subset of ${\bf\Sigma}_{\xi}\cap {\bf\Pi}_{\xi}\cap {\bf\Delta}$. So 
${\bf\Sigma}_{\xi}$, ${\bf\Pi}_{\xi}$ and ${\bf\Delta}$ are co-nowhere 
dense, and it remains to see that they are not open. It is enough to 
notice that $\emptyset$ is not in their interior. Look at the 
proof of Theorem 12; it shows that for each integer $m$, there is a 
subset $A_{m}$ of $\{s\in 5^{<\omega}/\vert s\vert \geq m\}$ such that 
$A_{m}^\infty\notin\borel$. But the argument in the proof of 
Proposition 11 shows that we can have the same thing in $n^{<\omega}$ 
for each $n\geq 2$. This gives the result because the sequence $(A_{m})$ tends 
to $\emptyset$.$\hfill\square$\bigskip

 We can say a bit more about ${\bf\Pi}_{1}$ and ${\bf\Sigma}_{2}$: 
 
\begin{prop} ${\bf\Pi}_{1}$, ${\it\Pi}_{1}$ and ${\bf\Sigma}_{2}$ are $\boratwo$-hard (so they 
are not $\bormtwo$).\end{prop} 

\bf\noindent Proof.\rm ~Consider the map $\phi$ defined in 
the proof of Proposition 5. By Proposition 2, if $\gamma\in P_{f}$, then 
$\phi (\gamma )^\infty$ is $\bormone$. Moreover, as $\phi (\gamma )$ is an antichain for the extension ordering, the decomposition into words of $\phi (\gamma )$ is unique. This shows 
that $\phi (\gamma )^\infty$ is $\Borel$, because 
$$\alpha\!\in\!\phi (\gamma )^\infty\Leftrightarrow~\exists 
\beta\!\in\!\Borel (\alpha )~[\big(\forall m\! >\! 0~~\beta (m)\! >\! 
0\big)~\mbox{and}~\big(\forall q\!\in\!\omega~~\pi (\alpha ,\beta ,q)\!\in\!\phi (\gamma )\big)].$$
So $\phi (\gamma )\in {\it\Pi}_{1}$ if $\gamma\in P_{f}$. So the preimage of any 
of the sets in the statement  by $\phi$ is $P_{f}$, and the result follows.$\hfill\square$\bigskip\smallskip

\noindent\bf {\Large 4 Which sets are $\omega$-powers?}\bigskip\rm

 Now we come to Question (3). Let us specify what we mean by ``codes 
for $\Gamma$-sets", where $\Gamma$ is a given class, and fix some notation.\bigskip

\noindent $\bullet$ For the Borel classes, we will essentially consider the 
$2^\omega$-universal sets used in [K1] (see Theorem 22.3). 
For $\xi\geq 1$, ${\cal U}^{\xi, {\cal A}}$ (resp. ${\cal U}^{\xi, {\cal 
M}}$) is $2^\omega$-universal for $\boraxi (n^\omega )$ (resp. $\bormxi 
(n^\omega )$). So we have\bigskip

\noindent  - ${\cal U}^{1, {\cal A}}=\{(\gamma ,\alpha )\in 2^\omega\!\times\! 
n^\omega /\exists p\in\omega ~\gamma (p)=0~\mbox{and}~s^n_{p}\prec\alpha\}$, where $(s^n_{p})_{p}$ enumerates $n^{<\omega}$.\smallskip

\noindent - ${\cal U}^{\xi, {\cal M}}=\neg ~{\cal U}^{\xi, {\cal A}}$, for each $\xi\geq 1$.\smallskip

\noindent - ${\cal U}^{\xi, {\cal A}}=\{(\gamma ,\alpha )\in 2^\omega\!\times\! 
n^\omega / \exists p\in\omega ~\big((\gamma )_{p},\alpha\big)\in 
{\cal U}^{\eta, {\cal M}}\}$ if $\xi =\eta +1$.\smallskip

\noindent - ${\cal U}^{\xi, {\cal A}}=\{(\gamma ,\alpha )\in 2^\omega\!\times\! 
n^\omega / \exists p\in\omega ~\big((\gamma )_{p},\alpha\big)\in 
{\cal U}^{\eta_{p}, {\cal M}}\}$ if $\xi$ is the limit of the strictly 
increasing sequence of odd ordinals $(\eta_{p})$.\bigskip

\noindent $\bullet$ For the class $\ana$, we fix some bijection $p\mapsto 
\big((p)_{0},(p)_{1}\big)$ between $\omega$ and $\omega^2$. We set
$$(\gamma ,\alpha )\!\in\! {\cal U}~\Leftrightarrow ~
\exists\beta\!\in\! 2^\omega ~\big(\forall m~\exists p\!\geq\! m~~\beta (p)\! 
=\! 1\big)~\mbox{and}~\big(\forall p~~[\gamma (p)\! =\! 1~\mbox{or}~s^2_{(p)_{0}}
\!\not\prec\!\beta ~\mbox{or}~s^n_{(p)_{1}}\!\not\prec\!\alpha ]\big).$$
It is not hard to see that $\cal U$ is $2^\omega$-universal for 
$\ana (n^\omega )$, and we use it here because of the compactness of 
$2^\omega\!\times\! n^\omega$, rather than the $\omega^\omega$-universal set for 
$\ana (n^\omega )$ given in [K1] (see Theorem 14.2).

\vfill\eject

\noindent $\bullet$ For the class $\borel$, it is different because there is no 
universal set. But we can use the $\Ca$ set of codes 
$D\subseteq 2^\omega$ for the Borel sets in [K1] (see Theorem 35.5). We may 
assume that $D$, $S$ and $P$ are effective, by [M].\bigskip

\noindent $\bullet$ The sets we are interested in are the following:
$${\cal A}_{\xi}:=\{\gamma\!\in\! 2^\omega /{\cal U}^{\xi, {\cal 
A}}_{\gamma}~\mbox{is~an}~\omega\mbox{-power}\}\mbox{,}~~~{\cal M}_{\xi}:=
\{\gamma\!\in\! 2^\omega /{\cal U}^{\xi, {\cal M}}_{\gamma}~\mbox{is~an}~\omega\mbox{-power}\}$$
$${\cal B}:=\{ d\!\in\! D /D_{d}~\mbox{is~an}~\omega\mbox{-power}\}\mbox{,}$$
$${\cal A}:=\{\gamma\!\in\! 2^\omega /{\cal U}_{\gamma}~\mbox{is~an}~\omega\mbox{-power}\}.$$
As we mentionned in the introduction, Lemma 13 is also related to 
Question (3). A rough answer to this question is ${{\it\Sigma}^1_{3}}$. 
Indeed, we have, for $\gamma\in 2^\omega$,
$$\gamma\!\in\! {\cal A}~\Leftrightarrow ~\exists A\!\in\! 
2^{n^{<\omega}}~\forall\alpha\!\in\! n^\omega ~\big([(\gamma ,\alpha )\!\notin\! 
{\cal U}~\mbox{or}~\alpha\!\in\! A^\infty ]~\mbox{and}~[\alpha\!\notin\! 
A^\infty ~\mbox{or}~(\gamma ,\alpha )\!\in\! {\cal U}]\big).$$
With Lemma 13, we have a better estimation of the complexity of $\cal B$: it 
is ${\it\Sigma}^1_2$. Indeed, for $d\in D$,\bigskip

\leftline{$D_{d}~\mbox{is~an}~\omega\mbox{-power}~\Leftrightarrow ~
\exists A\!\in\! 2^{n^{<\omega}}~\forall\alpha\!\in\! n^\omega ~
\big(\big[(d ,\alpha )\!\notin\! S~~\mbox{or}~~\exists\beta\!\in\! 
\Borel (A,d ,\alpha )$}\bigskip

\rightline{$[\big(\forall m\! >\! 0~~\beta (m)\! >\! 0\big)~\mbox{and}~
\big(\forall q\!\in\!\omega~~\pi (\alpha ,\beta ,q)\!\in\! A\big)]\big]
~\mbox{and}~[\alpha\!\notin\! A^\infty ~\mbox{or}~(d ,\alpha )\!\in\! 
P]\big).$}\bigskip

\noindent This argument also shows that ${\cal A}_{\xi}$ and ${\cal M}_{\xi}$ 
are ${\bf\Sigma}^1_2$. We can say more about these two sets.

\begin{prop} If $1\leq\xi<\omega_1$, then 
${\cal A}_{\xi}$ and ${\cal M}_{\xi}$ are ${\bf\Sigma}^1_2\setminus D_{2}(\boraone )$ co-meager subsets of $2^\omega$. If moreover $\xi =1$, then they are co-nowhere dense.\end{prop} 

\bf\noindent Proof.\rm ~We set $E_{1}:=
\{\gamma\!\in\! 2^\omega /{\cal U}^{1, {\cal A}}_{\gamma}\! =\! n^\omega\}$, 
$E_{\eta +1}:=\{\gamma\!\in\! 2^\omega /\forall p~(\gamma )_{p}\!\in\! 
E_{\eta}\}$ if $\eta\geq 1$, and $E_{\xi}:=\{\gamma\!\in\! 2^\omega /
\forall p~(\gamma )_{p}\!\in\! E_{\eta_{p}}\}$ (where $(\eta_{p})$ is a strictly 
increasing sequence of odd ordinals cofinal in the limit ordinal $\xi$). 
If $s\in 2^{<\omega}$, then we set $\gamma (p)=s(p)$ if $p<\vert s\vert $, $0$ 
otherwise. Then $s\prec\gamma$ and ${\cal U}^{1, {\cal 
A}}_{\gamma}=n^\omega$, so $E_{1}$ is dense. If $\gamma_{0}\in E_{1}$, then for all 
$\alpha\in n^\omega$ we can find an integer $p$ such 
that $\gamma_{0}(p)=0$ and $s^n_{p}\prec\alpha$. By compactness of 
$n^\omega$ we can find a finite subset $F$ of $\{p\!\in\!\omega /\gamma_{0}(p)\! =\! 0\}$ such that for each $\alpha\in n^\omega$, $s^n_{p}\prec\alpha$ for some $p\in F$. Now 
$\{\gamma\!\in\! 2^\omega /\forall p\!\in\! F~\gamma (p)\! =\! 0\}$ is an open 
neighborhood of $\gamma_{0}$ and a subset of $E_{1}$. So $E_{1}$ is an open 
subset of $2^\omega$. Now the map $\gamma\mapsto (\gamma)_{p}$ is 
continuous and open, so $E_{\eta +1}$ and $E_{\xi}$ are dense $G_{\delta}$ 
subsets of $2^\omega$. Then we notice that $E_{\xi}$ is a subset of 
$\{\gamma\!\in\! 2^\omega /{\cal U}^{\xi, {\cal A}}_{\gamma}\! =\! n^\omega\}$ 
(resp., $\{\gamma\!\in\! 2^\omega /{\cal U}^{1, {\cal 
A}}_{\gamma}\! =\!\emptyset\}$) if $\xi$ is odd (resp., even). Indeed, this 
is clear for $\xi=1$. Then we use the formulas 
${\cal U}^{\eta +1, {\cal A}}_{\gamma} = \bigcup_{p}~\neg ~ 
{\cal U}^{\eta , {\cal A}}_{(\gamma )_{p}}$ and ${\cal U}^{\xi , 
{\cal A}}_{\gamma} = \bigcup_{p}~\neg ~{\cal U}^{\eta_{p} , 
{\cal A}}_{(\gamma )_{p}}$, and by induction we are done. As 
$\emptyset$ and $n^\omega$ are $\omega$-powers, we get the results 
about Baire category. Now it remains to see that ${\cal A}_{\xi}$ and 
${\cal M}_{\xi}$ are not open. But by induction 
again $1^\infty\in {\cal A}_{\xi}\cap {\cal M}_{\xi}$, so it is enough to see 
that $1^\infty$ is not in the interior of these sets.\bigskip

\noindent $\bullet$ Let us show that, for $O\in\borone(n^\omega )\setminus\{\emptyset ,n^\omega\}$ 
and for each integer $m$, we can find $\gamma$, $\gamma'\in 2^\omega$ 
such that $\gamma (j)=\gamma' (j)=1$ for $j<m$, 
${\cal U}^{\xi , {\cal A}}_{\gamma} = O$ and 
${\cal U}^{\xi , {\cal M}}_{\gamma'} = O$.\bigskip

 For $\xi =1$, write $O=\bigcup_{p}~N_{s^n_{q_{k}}}$, where $q_{k}\geq 
m$. Let $\gamma (q):=0$ if there exists $k$ such that $q=q_{k}$, $\gamma 
(q):=1$ otherwise. The same argument applied to $\check O$ gives the 
complete result for $\xi =1$.

\vfill\eject

 Now we argue by induction. Let 
$\gamma_{p}\in 2^\omega$ be such that $\gamma_{p}(q)=1$ for $<\! p,q\! ><m$ 
and ${\cal U}^{\eta , {\cal M}}_{(\gamma )_{p}} = O$. Then define 
$\gamma$ by $\gamma (<\! p,q\! >):=\gamma_{p}(q)$; we have $\gamma (j)=1$ 
if $j<m$ and ${\cal U}^{\eta +1, {\cal A}}_{\gamma} = \bigcup_{p}~
{\cal U}^{\eta , {\cal M}}_{(\gamma )_{p}} = O$. The argument with 
$\check O$ still works. The argument is similar for limit  ordinals.\bigskip

\noindent $\bullet$ Now we apply this fact to $O := N_{(0)}$. This gives 
$\gamma_{p},~\gamma'_{p}\in N_{1^p}$ such that 
${\cal U}^{\xi , {\cal A}}_{\gamma _{p}}=N_{(0)}$ and 
${\cal U}^{\xi , {\cal M}}_{\gamma' _{p}}=N_{(0)}$. 
But $(\gamma_{p})$, $(\gamma'_{p})$ tend to $1^\infty$, 
$\gamma_{p}\notin {\cal A}_{\xi}$ and $\gamma'_{p}\notin {\cal M}_{\xi}$.$\hfill\square$

\begin{cor} ${\cal A}_{1}$ is 
$\check D_{2}(\boraone)\setminus D_{2}(\boraone)$. In particular, ${\cal A}_{1}$ is 
$\check D_{2}(\boraone)$-complete.\end{cor}

\bf\noindent Proof.\rm ~By the preceding proof, it is 
enough to see that ${\cal A}_{1}\setminus\{1^\infty\}$ is open. So 
let ${\gamma_{0}\in {\cal A}_{1}\setminus\{1^\infty\}}$,  
$p_{0}$ in $\omega$ with $\gamma_{0}(p_{0})=0$, and 
$A_{0}\subseteq n^{<\omega}$ with ${{\cal U}^{1 , {\cal A}}_{\gamma _{0}}=A_{0}^\infty}$. If $\alpha\in 
n^\omega$, then ${s^n_{p_{0}}}\alpha\in {\cal U}^{1 , {\cal A}}_{\gamma 
_{0}}$, so we can find $m>0$ such that ${\alpha -\alpha\lceil 
m\in A_{0}^\infty}$; thus there exists an integer $p$ such that 
$\gamma_{0}(p)=0$ and ${s^n_{p}\prec \alpha -\alpha\lceil m}$. By 
compactness of $n^\omega$, there are finite sets 
$F\subseteq\omega\setminus\{ 0\}$ and 
${G\subseteq\{ p\!\in\!\omega /\gamma_{0}(p)\! =\! 0\}}$ 
such that ${n^\omega =\bigcup_{m\in F,p\in G}~\{\alpha\!\in\! n^\omega 
/s^n_{p}\!\prec\!\alpha\! -\!\alpha\lceil m\}}$. 

 We set ${A_{\gamma}:=\{s\in n^{<\omega}/\exists 
p~\gamma (p)=0~\mbox{and}~s^n_{p}\prec s\}}$ for $\gamma\in 2^\omega$, so that 
$A_{\gamma}^\infty\subseteq {\cal U}^{1 , {\cal A}}_{\gamma}$. Assume 
that $\gamma (p)=0$ for each $p\in G$ and let $\alpha\in 
{\cal U}^{1 , {\cal A}}_{\gamma}$. Let $p^0\in\omega$ be such that 
$\gamma (p^0)=0$ and $s^n_{p^0}\prec\alpha$. We can find  
$m_{0}>0$ and $p^1\in G$ such that ${s^n_{p^1}\prec\alpha 
-\alpha\lceil (\vert s^n_{p^0}\vert +m_{0})}$, and 
$\alpha\lceil (\vert s^n_{p^0}\vert +m_{0})\in A_{\gamma}$. Then we can find  
$m_{1}>0$ and $p^2\in G$ such that ${s^n_{p^2}\prec\alpha 
-\alpha\lceil (\vert s^n_{p^0}\vert +m_{0}+\vert s^n_{p^1}\vert +m_{1})}$, and 
$$\alpha\lceil (\vert s^n_{p^0}\vert +m_{0}+\vert s^n_{p^1}\vert +m_{1})-\alpha\lceil 
(\vert s^n_{p^0}\vert +m_{0})\in A_{\gamma}.$$ 
And so on. Thus $\alpha\in A_{\gamma}^\infty$ and 
$\{\gamma\!\in\! 2^\omega /\forall p\!\in\! G~\gamma (p)\! =\! 0\}$ is a clopen neighborhood of $\gamma_{0}$ and a subset of ${\cal A}_{1}$.$\hfill\square$

\begin{prop} ${\cal A}$ is 
${{\it\Sigma}^1_{3}}\setminus D_{2}(\boraone)$ and is co-nowhere dense.\end{prop} 

\bf\noindent Proof.\rm ~Let $U:=\{\gamma\!\in\! 
2^\omega /\forall\beta\!\in\! 2^\omega ~\forall\alpha\!\in\! n^\omega ~\exists 
p~~[\gamma (p)\! =\! 0~\mbox{and}~s^2_{(p)_{0}}\!\prec\!\beta 
~\mbox{and}~s^n_{(p)_{1}}\!\prec\!\alpha ]\}$. By compactness of 
$2^\omega\times n^\omega$, $U$ is a dense open subset of $2^\omega$. 
Moreover, if $\gamma\in U$, then ${\cal U}_{\gamma}=\emptyset$, so 
$U\subseteq {\cal A}$ and $\cal A$ is co-nowhere dense. It remains to 
see that $\cal A$ is not open, as in the proof of Proposition 18. 
As ${\cal U}_{1^\infty}=n^\omega$, $1^\infty\in {\cal A}$. Let $p$ be 
an integer satisfying $s^2_{(p)_{0}}=\emptyset$ and $s^n_{(p)_{1}}=0^q$. We 
set $\gamma_{p}(m):=0$ if and only if $m=p$, and also $P_{\infty} :=
\{\alpha\!\in\! 2^\omega /\forall r~\exists m\!\geq\! r~\alpha (m)\! =\! 1\}$. 
Then $(\gamma_{p})$ tends to $1^\infty$ and we have 
$$\begin{array}{ll}
{\cal U}_{\gamma_{p}}\!\!\!
& =\{\alpha\!\in\! n^\omega /\exists\beta\!\in\! P_{\infty}~\forall m~~m\!\not=\! p~\mbox{or}~s^2_{(m)_{0}}\!\not\prec\!\beta ~\mbox{or}~s^n_{(m)_{1}}\!\not\prec\!\alpha\} \cr 
& =\{\alpha\!\in\! n^\omega /\exists\beta\!\in\! P_{\infty}~(\beta ,\alpha 
)\!\notin\! 2^\omega\!\times\! N_{0^q}\} = \neg ~N_{0^q}. 
\end{array}$$
So $\gamma_{p}\notin {\cal A}$.$\hfill\square$\bigskip\smallskip

\noindent\bf {\Large 5 Ordinal ranks and $\omega$-powers.}\rm\bigskip

\noindent\bf Notation.\rm ~The fact that the $\omega$-powers are $\ana$ 
implies the existence of a co-analytic rank on the complement of $A^\infty$ 
(see 34.4 in [K1]). We will consider a natural one,
defined as follows. We set, for $\alpha\in n^\omega$, 
$T_A(\alpha):=\{S\!\in\! (A^{-})^{<\omega}/S^{*}\!\prec\!\alpha\}$. This is a tree on 
$A^-$, which is well founded if and only if 
$\alpha\notin A^\infty$.

\vfill\eject

 The rank of this tree is the announced rank 
$R_A:\neg ~A^\infty\rightarrow\omega_1$ (see page 10 in [K1]): 
we have $R_A(\alpha ):=\rho\big(T_A(\alpha )\big)$. Let 
$\phi :A^{-}\rightarrow\omega$ be one-to-one, and 
$\tilde\phi (S):=\big(\phi [S(0)],\ldots ,\phi [S(\vert s\vert -1)]\big)$ for 
$S\in (A^{-})^{<\omega}$. This allows us to define the map $\Phi$ from the set 
of trees on $A^{-}$ into the set of trees on $\omega$, which associates 
$\{\tilde\phi (S)/S\!\in\! T\}$ to $T$. As $\tilde\phi$ is one-to-one, 
$\Phi$ is continuous: 
$$t\in\Phi (T)\Leftrightarrow t\in\tilde\phi [(A^{-})^{<\omega}]~\mbox{and}~\tilde\phi^{-1}(t)\in T.$$ 
Moreover, $T$ is well-founded if and only if $\Phi (T)$ is well-founded. Thus, if 
$\alpha\notin A^\infty$, then we have $\rho \big(T_A(\alpha )\big)=\rho 
\big(\Phi [T_A(\alpha )]\big)$ because $\tilde\phi$ is strictly monotone 
(see page 10 in [K1]). Thus 
$R_A$ is a co-analytic rank because the function from $n^\omega$ into the set 
of trees on $\omega^{<\omega}$ which associates $\Phi [T_A(\alpha )]$ to $\alpha$ 
is continuous, and because the rank of the well-founded trees on $\omega$ 
defines a co-analytic rank (see 34.6 in [K1]). We set 
$$R(A) := \mbox{sup}\{R_A(\alpha )/\alpha\!\notin\! A^\infty\}.$$
By the boundedness theorem, $A^\infty$ is Borel if and only if 
$R(A)<\omega_1$ (see 34.5 and 35.23 in [K1]). 
We can ask the question of the link between the complexity of $A^\infty$ 
and the ordinal $R(A)$ when $A^\infty$ is Borel.

\begin{prop} If $\xi<\omega_1$, $r\in\omega$ 
and $R(A)=\omega.\xi +r$, then $A^\infty\in {\bf\Sigma}^0_{2.\xi +1}$.\end{prop} 

\bf\noindent Proof.\rm ~The reader should see [L] for operations on ordinals.\bigskip

\noindent $\bullet$ If $0<\lambda<\omega_1$ is a limit ordinal, then let $(\lambda_q)$ be a strictly 
increasing co-final sequence in $\lambda$, with $\lambda_q = \omega .\theta
+q$ if $\lambda = \omega .(\theta +1)$, and $\lambda_q = \omega .\xi_q$ 
if $\lambda = \omega .\xi$, where $(\xi_q)$ is a strictly increasing 
co-final sequence in the limit ordinal $\xi$ otherwise. By induction, we define 
$$\begin{array}{ll} 
E_0\!\!\!\! 
& := \{\alpha\!\in\! n^\omega/\forall s\!\in\! A^{-}~s\!\not\prec\!\alpha\}\mbox{,}\cr 
E_{\theta +1}\!\!\!\! 
& := \{\alpha\!\in\! n^\omega/\forall s\!\in\! A^{-}~s\!\not\prec\!
\alpha~~\mbox{or}~~\alpha \! -\! s\!\in\! E_\theta\}\mbox{,}\cr
E_\lambda\!\!\!\! 
& :=\{\alpha\!\in\! n^\omega/\forall s\!\in\! A^{-}~s\!\not\prec\!
\alpha~~\mbox{or}~~\exists q\!\in\!\omega~~\alpha\! -\! s\!\in\! E_{\lambda_q}\}. 
\end{array}$$
$\bullet$ Let us show that $E_{\omega.\xi +r}\in {\bf\Pi}^0_{2.\xi +1}$. 
We may assume that $\xi\not=0$ and that $r=0$. If $\xi =\theta +1$, then 
$E_{\lambda _q}\in {\bf\Pi}^0_{2.\theta +1}$ by induction hypothesis, thus 
$E_{\omega.\xi +r}\in {\bf\Pi}^0_{2.\theta +3}={\bf\Pi}^0_{2.\xi +1}$. Otherwise, 
$E_{\lambda _q}\in {\bf\Pi}^0_{2.\xi_q +1}$ by induction hypothesis, thus 
$E_{\omega.\xi +r}\in {\bf\Pi}^0_{\xi +1}={\bf\Pi}^0_{2.\xi +1}$.\bigskip

\noindent $\bullet$ Let us show that if $\alpha\in A^\infty$, then 
$\alpha\notin E_{\omega.\xi +r}$. If $\xi = r = 0$, it is clear. If 
$r=m+1$ and $s\in A^{-}$ satisfies $s\prec\alpha$ and $\alpha -
s\in A^\infty$, then we have $\alpha -s\notin E_{\omega.\xi +m}$ by induction 
hypothesis, thus $\alpha\notin  E_{\omega.\xi +r}$. If $r=0$ and $s\in A^{-}$ 
satisfies $s\prec\alpha$ and $\alpha - s\in
A^\infty$, then we have $\alpha -s\notin E_{\lambda_q}$ for each integer 
$q$, by induction hypothesis, thus $\alpha\notin E_{\omega.\xi +r}$.\bigskip

\noindent $\bullet$ Let $s\in A^{-}$ such that $s\prec\alpha\notin A^\infty$. We have 
$$\begin{array}{ll} 
\rho \big(T_A(\alpha -s)\big)\!\!\!\! 
& = \mbox{sup}\{\rho_{T_A(\alpha -s)}(t)+1~/~t\in T_A(\alpha -s)\}\cr 
& \leq\mbox{sup}\{\rho_{T_A(\alpha )}\big( (s) t\big)+1~/~(s) t\in T_A(\alpha )\}
\cr & \leq \rho_{T_A(\alpha )}\big( (s)\big) +1\cr 
& \leq \rho_{T_A(\alpha )}(\emptyset ) < \rho\big(T_A(\alpha )\big).
\end{array}$$

\vfill\eject

 The first inequality comes from the fact that the map from $T_{A}(\alpha -s)$ into $T_{A}(\alpha )$, which associates $(s)t$ to $t$ is strictly monotone (see page 10 in [K1]). We have  
$${\rho\big(T_A(\alpha )\big)\geq [\mbox{sup}\{\rho\big(T_A(\alpha -s)\big)~/~s\in A^{-},~s\prec\alpha\}]
+1}.$$
Let us show that we actually have equality. We have 
$$\rho\big(T_A(\alpha )\big)=\rho_{T_A(\alpha )}(\emptyset )+1=\mbox{sup}\{\rho_{T_A(\alpha )}\big( (s)\big)+1~/~s\in A^{-},~s\prec\alpha\}+1.$$
Therefore, it is enough to notice that if $s\! \in\!  A^-$ and $s\! \prec\!\alpha$, then 
${\rho_{T_A(\alpha )}\big( (s)\big)\! \leq\!  \rho_{T_A(\alpha -s)}(\emptyset )}$. But this comes from the fact that the map from ${\{S\in T_{A}(\alpha )~/~S(0)=s\}}$ into $T_A(\alpha -s)$, which 
associates $S-(s)$ to $S$, preserves the extension ordering (see page 352 in [K1]).\bigskip

\noindent $\bullet$ Let us show that, if $\alpha\notin A^\infty$, then ``$\rho \big(T_A(\alpha )
\big)\leq \omega.\xi +r+1$" is equivalent to ``$\alpha\in E_{\omega.\xi +r}$". 
We do it by induction on $\omega.\xi +r$. If $\xi=r=0$, then it is clear.
If $r=m+1$, then ``$\rho \big(T_A(\alpha )\big)\leq \omega.\xi +r+1$" is 
equivalent to ``$\forall s\in A^{-}$,
$s\not\prec\alpha$ or $\rho \big(T_A(\alpha -s)\big)\leq \omega.\xi +m+1$", 
by the preceding point. This is equivalent to ``$\forall s\in A^{-}$, 
$s\not\prec\alpha$ or $\alpha -s\in E_{\omega.\xi +m}$", which is 
equivalent to ``$\alpha\in E_{\omega.\xi +r}$". If $r=0$, then ``$\rho
\big(T_A(\alpha )\big)\leq \omega.\xi +r+1$" is equivalent to 
``$\forall s\in A^{-}$, $s\not\prec\alpha$ or there exists an integer $q$ 
such that $\rho\big(T_A(\alpha -s)\big)\leq\lambda_q+1$". This is equivalent to 
``$\forall s\in A^{-}$, $s\not\prec\alpha$ or there exists an integer $q$ 
such that $\alpha -s\in E_{\lambda_q}$", which is equivalent to 
``$\alpha\in E_{\omega.\xi +r}$".\bigskip

\noindent $\bullet$ If $\alpha\notin A^\infty$, then $\rho\big(T_A(\alpha )\big)\leq 
\omega.\xi +r+1$. By the preceding point, $\alpha\in
E_{\omega.\xi +r}$. Thus we have $A^\infty = \neg ~ E_{\omega.\xi +r}\in
{\bf\Sigma}^0_{2.\xi +1}$.$\hfill\square$\bigskip
 
  We can find an upper bound for the rank $R$, for some Borel classes:
 
\begin{prop} (a) $A^\infty=n^\omega$ 
if and only if $R(A)=0$.\smallskip

\noindent (b) If $A^\infty =\emptyset$, then $R(A)=1$.\smallskip

\noindent (c) If $A^\infty\in\borone$, then $R(A)<\omega$, and there exists 
$A_{p}\subseteq 2^{<\omega}$ such that $A_{p}^\infty\in\borone$ 
and ${R(A_{p})=p}$ for each integer $p$.\smallskip

\noindent (d) If $A^\infty\in\bormone$, then $R(A)\leq\omega$, and  
$(A^\infty\notin\boraone\Leftrightarrow R(A)=\omega)$.\end{prop}
 
\bf\noindent Proof.\rm ~(a) If $\alpha\notin A^\infty$, 
then $\emptyset\in T_{A}(\alpha )$ and ${\rho\big(T_{A}(\alpha )\big)\!\geq\!
\rho_{T_{A}(\alpha )}(\emptyset )\! +\! 1\!\geq\! 1}$.\bigskip

\noindent (b) We have $T\! _{\! A}(\alpha )\! =\!\{\emptyset\}$ for each $\alpha$, and 
${\rho\big(T_{A}(\alpha )\big)\! =\! \rho_{T\! _{\! A}(\alpha )}(\emptyset)\! +\! 1\! =\! 1}$.\bigskip

\noindent (c) By compactness, there exists $s_{1},\ldots,s_{p}\in n^{<\omega}$ 
such that $A^\infty=\bigcup_{1\leq m\leq p}~N_{s_{m}}\in\borone$. If  
$\alpha\notin A^\infty$, then we have $N_{\alpha\lceil\mbox{max}_{1\leq m\leq 
p}~\vert s_{m}\vert }\subseteq\neg ~A^\infty$, thus $\rho\big(T_{A}(\alpha )\big)\leq 
\mbox{max}_{1\leq m\leq p}~\vert s_{m}\vert +1<\omega$. So we get the first point. 
To see the second one, we set $A_{0} := 2^{<\omega}$. If $p>0$, then we set 
$$A_{p}:= \{ 0^2\}\cup\bigcup_{q\leq p}~
\{s\!\in\! 2^{<\omega}/0^{2q}1\!\prec\! s\}\cup 
\{s\!\in\! 2^{<\omega}/0^{2p+1}\!\prec\! s\}.$$
Then $A_{p}^\infty =\bigcup_{q\leq p}~N_{0^{2q}1}\cup 
N_{0^{2p+1}}\in\borone$. If $\alpha_{p} := 0^{2p-1}1^\infty$, then  
$\rho\big(T_{A_{p}}(\alpha_{p} )\big)=p$. If $\alpha\notin A_{p}^\infty$, then 
$\rho\big(T_{A_{p}}(\alpha )\big)\leq p$.

\vfill\eject

\noindent (d) If $A^\infty\in\bormone$ and $\alpha\notin A^\infty$, then let $s\in 
n^{<\omega}$ with $\alpha\in N_{s}\subseteq\neg ~A^\infty$. 
Then ${\rho\big(T_{A}(\alpha )\big)\leq \vert s\vert +1}$. Thus $R(A)\leq\omega$. 
If $A^\infty\notin\boraone$, then we have $R(A)\geq\omega$, by Proposition 21. 
Thus $R(A)=\omega$. Conversely, we apply (c).$\hfill\square$\bigskip

\noindent\bf Remark.\rm ~Notice that it is not true that if the Wadge class $<A^\infty >$, having 
$A^\infty$ as a complete set, is a subclass of $<B^\infty >$, then 
$R(A)\leq R(B)$. Indeed, for $A$ we take the example $A_{2}$ in (c), and for $B$ we take the example for 
$\boraone$ that we met in the proof of Proposition 11. If we exchange the roles of $A$ and $B$, then we see that the converse is also false. This example $A$ for $\boraone$ shows that 
Proposition 21 is optimal for $\xi=0$ since $R(A)=1$ and $A^\infty\in\boraone\setminus\bormone$. We can say more: it is not true that if 
$A^\infty =B^\infty$, then $R(A)\leq R(B)$. We use again (c): 
we take $A:=A_{2}$ and $B:=A\setminus\{ 0^2\}$. We have $A^\infty =B^\infty = 
A_{2}^\infty$, $R(A)=2$ and $R(B)=1$.

\begin{prop} For each $\xi <\omega_1$, 
there exists $A_{\xi}\subseteq 2^{<\omega}$ with  
$A_{\xi}^{\infty}\in\boraone$ and $R(A_{\xi})\geq\xi$.\end{prop}

\bf\noindent Proof.\rm ~We use the notation in the 
proof of Theorem 15. Let $T\!\in\! {\cal T}$, and ${\varphi : T\rightarrow 
T_{\Phi'(T)}(\alpha_{0})}$ defined by the formula $\varphi (s):=\big(\phi (s\lceil 0),
\ldots ,\phi (s\lceil\vert s\vert -1)\big)$. Then $\varphi$ is strictly monotone. If 
$T\in WF$, then $\alpha_{0}\notin\big(\Phi'(T)\big)^\infty$ and 
$T_{\Phi'(T)}(\alpha_{0})\in WF$. In this case, 
${\rho\big(T\big)\leq\rho \big(T_{\Phi'(T)}(\alpha_{0})\big)=
R_{\Phi'(T)}(\alpha_{0})}$ (see page 10 in [K1]). Let $T_{\xi}\in 
WF$ be a tree with rank at least $\xi$ (see 34.5 and 34.6 in [K1]). We set $A_{\xi}:=\Phi'(T_{\xi})$. 
It is clear that $A_{\xi}$ is what we were looking for.$\hfill\square$\bigskip
 
\noindent\bf Remark.\rm ~Let $\psi : 2^{n^{<\omega}}\!\rightarrow\!\{\mbox{Trees~on}~n^{<\omega}\}$ defined by $\psi (A)\! :=\! T_{A}(\alpha_{0})$, and $r :\neg ~I_{\alpha_{0}}\!\rightarrow\!\omega_{1}$ 
defined by $r(A):=\rho\big(T_{A}(\alpha_{0})\big)$. Then $\psi$ is continuous, thus $r$ is a 
$\ca$-rank on 
$${\psi^{-1}(\{\mbox{Well-founded~trees~on}~n^{<\omega }\})\! =\! 
\neg ~I_{\alpha_{0}}}.$$ 
By the boundedness theorem, the rank $r$ and $R$ are not bounded on $\neg ~I_{\alpha_{0}}$. Proposition 23 specifies this result. It shows that $R$ is not bounded on ${\it\Sigma}_{1}\setminus I_{\alpha_{0}}$.\bigskip\smallskip

\noindent\bf {\Large 6 The extension ordering.}\rm

\begin{prop} We equip $A$ with the 
extension ordering.\smallskip

\noindent (a) If $A\subseteq n^{<\omega}$ is an antichain, then $A^\infty$ is in 
$\{\emptyset\}\cup\{n^\omega\}\cup [\bormone\setminus\boraone]\cup 
[\Bormtwo(A)\setminus\boratwo]$, and any of these cases is possible.\smallskip

\noindent (b) If $A\subseteq n^{<\omega}$ has finite antichains, then $A^\infty\in\bormtwo$ (and is not 
$\boratwo$ in general).\end{prop}
  
\bf\noindent Proof.\rm ~Let $G:=\{\alpha\in n^\omega~/~\forall r~
\exists m~\exists p\geq r~~\alpha\lceil m\in [(A^{-})^p]^{*}\}$. 
Then $G\in\Bormtwo(A)$ and contains $A^\infty$. Conversely, if $\alpha\in G$, then we 
have $T_{A}(\alpha )\cap (A^{-})^p\not= \emptyset$ for each integer $p$, thus 
$T_{A}(\alpha )$ is infinite.\bigskip

\noindent (a) If $A$ is an antichain, then each sequence in $T_{A}(\alpha )$ has at most 
one extension in this tree adding one to the length. Thus $T_{A}(\alpha )$ is 
finite splitting. This implies that $T_{A}(\alpha )$ has an infinite branch if 
$\alpha\in G$, by K\" onig's lemma. Therefore $A^\infty =G\in\Bormtwo(A)$.

\vfill\eject

\noindent - If we take $A:=\emptyset$, then $A$ is an antichain and $A^\infty=\emptyset$.\smallskip

\noindent - If we take $A:=\{ (0),\ldots ,(n-1)\}$, then $A$ is an antichain and $A^\infty=n^\omega$.\smallskip

\noindent - If $A^\infty\notin\{\emptyset , n^\omega\}$, then $A^\infty\notin\boraone$. Indeed, let $\alpha_{0}\notin A^\infty$ and $s_{0}\in A^{-}$. By uniqueness 
of the decomposition into words of $A^{-}$, the sequence 
$({s_{0}^n}\alpha_{0})_{n}\subseteq n^\omega\setminus A^\infty$ 
tends to $s_{0}^\infty\in A^\infty$.\smallskip

\noindent - If we take $A:=\{ (0)\}$, then $A$ is an antichain and $A^\infty\! =\! \{0^\infty\}
\in\bormone\setminus\boraone$.\smallskip

\noindent - If $A$ is finite, then $A^\infty$ is $\bormone\setminus\boraone$ or is in 
$\{\emptyset ,n^\omega\}$, by the facts above and Proposition 2.\smallskip

\noindent - If $A$ is infinite, then $A^\infty\notin\boratwo$ because the map $c$ in 
the proof of Proposition 2 is an homeomorphism and $(A^{-})^\omega$ 
is not $K_{\sigma}$.\smallskip

\noindent - If $A:=\{0^k1/k\!\in\!\omega\}$, then $A$ is an antichain and 
${A^\infty=P_{\infty}}$, which is $\Bormtwo\setminus\boratwo$.\smallskip

\noindent (b) The intersection of $P_{\infty}$ with $N_{1}$ can be 
made with the chain $\{10^k/k\!\in\!\omega\}$. So let us assume that $A$ 
has finite antichains.\bigskip

\noindent $\bullet$ Let us show that $A$ is the union of a finite set and of a 
finite union of infinite subsets of sets of the form 
$A_{\alpha_m}:=\{s\in n^{<\omega}/s\!\prec\!\alpha_{m}\}$. Let us 
enumerate $A:=\{s_{r}/r\!\in\!\omega\}$. We construct a sequence $(A_{m})$, finite 
or not, of subsets of $A$. We do it by induction on $r$, to decide in 
which set $A_{m}$ the sequence $s_r$ is. First, $s_0\in A_{0}$. Assume that $s_{0}, \ldots 
,s_{r}$ have been put into $A_{0},\ldots ,A_{p_{r}}$, with $p_{r}\leq r$ 
and $A_{m}\cap \{s_{0}, \ldots ,s_{r}\}\not=\emptyset$ if $m\leq p_{r}$. 
 We choose $m\leq p_{r}$ minimal such that $s_{r+1}$ is compatible 
with all the sequences in $A_{m}\cap \{s_{0}, \ldots ,s_{r}\}$, we put $s_{r+1}$ 
into $A_{m}$ and we set $p_{r+1}:=p_{r}$ if possible. Otherwise, 
we put $s_{r+1}$ into $A_{p_{r}+1}$ and we set $p_{r+1}:=p_{r}+1$.\bigskip

 Let us show that there are only finitely many infinite $A_{m}$'s. If $A_{m}$ is 
infinite, then there exists a unique sequence $\alpha_{m}\in n^\omega$ such 
that $A_{m}\subseteq A_{\alpha_m}$. Let us argue by contradiction: 
there exists an infinite sequence $(m_{q})_{q}$ such that $A_{m_{q}}$ is 
infinite. Let $t_{0}$ be the common beginning of the $\alpha_{m_{q}}$'s. 
There exists $\varepsilon_{0}\in n$ such that 
$N_{t_{0}\varepsilon_{0}}\cap\{\alpha_{m_{q}}/q\!\in\!\omega\}$ is 
infinite. We choose a sequence $u_{0}$ in $A$ extending $t_{0}
\mu_{0}$, where $\mu_{0}\not=\varepsilon_{0}$. Then we do it again: let 
$t_{0}\varepsilon_{0}t_{1}$ 
be the common beginning of the elements of $N_{t_{0}\varepsilon_{0}}\cap
\{\alpha_{m_{q}}/q\!\in\!\omega\}$. There exists $\varepsilon_{1}\in n$ such 
that $N_{t_{0}\varepsilon_{0} t_{1}\varepsilon_{1}}\cap
\{\alpha_{m_{q}}/q\!\in\!\omega\}$ is infinite. We choose a sequence 
$u_{1}$ in $A$ extending $t_{0}\varepsilon_{0} t_{1} 
\mu_{1}$, where $\mu_{1}\not=\varepsilon_{1}$. The sequence $(u_{l})$ is an 
infinite antichain in $A$. But this is absurd. Now let us choose the 
longest sequence in each nonempty finite $A_{m}$; this gives an 
antichain in $A$ and the result.\bigskip

\noindent $\bullet$ Now let $\alpha\in G$. There are two cases. Either for each 
$m$ and for each integer $k$, $\alpha\lceil k\notin [A^{<\omega}]^{*}$ or 
$\alpha -\alpha\lceil k\not= \alpha_{m}$. In this case, $T_{A}(\alpha )$ 
is finite splitting. As $T_{A}(\alpha )$ is infinite, $T_{A}(\alpha )$ has an 
infinite branch witnessing that $\alpha\in A^\infty$, by K\" onig's lemma. 
Otherwise, ${\alpha\in \bigcup_{s\in [A^{<\omega}]^{*},m} \{s\alpha_{m}\}}$, 
which is countable. Thus $G\setminus A^\infty\in\boratwo$ and 
${A^\infty =G\setminus (G\setminus A^\infty)\!\in\!\bormtwo}$.$\hfill\square$\bigskip\smallskip

\noindent\bf {\Large 7 Examples.}\rm\bigskip

\noindent $\bullet$ We have seen examples of subsets 
$A$ of $2^{<\omega}$ such that $A^\infty$ is complete for the classes 
$\{\emptyset\}$, $\{n^\omega\}$, $\borone$, $\boraone$, $\bormone$, $\bormtwo$ and $\ana$. We will give some more examples, for some classes of Borel sets. Notice that to show that a set in such a non self-dual class is complete, it is enough to show that it is true (see 21.E, 22.10 and 22.26 in [K1]).

\vfill\eject

\noindent $\bullet$ For the class $\boraone\oplus\bormone:=\{(U\cap O)\cup (F\setminus 
 O)~/~U\in\boraone,~ O\in\borone,~ F\in\bormone\}$, we can take 
${A\! :=\! \{s\! \in\!  2^{<\omega}\! /0^21\! \prec\!  s~\mbox{or}~s\! =\! 
0^2~\mbox{or}~\exists p\! \in\! 
\omega~10^p1\! \prec\!  s\}}$, since ${A^\infty\!  =\! \{0^\infty\}\! \cup\!  
\bigcup_{q}\!  N_{0^{2q+2}1}\! \cup\!  N_{1}\!\!  \setminus\!\!  \{ 10^\infty\}}$.\bigskip

\noindent $\bullet$ For the class $\check D_{2}(\boraone):=\{U\cup F~/~U\in\boraone ,~ F\in\bormone\}$, we can take Example 9 in [St2]: 
$A := \{s\in 2^{<\omega}~/~0\prec s~\mbox{or}~\exists ~q\in\omega 
~~{(101)^q}1^3\prec s~\mbox{or}~s=10^2\}$. We have 
$$A^\infty = 
\bigcup_{p\in\omega}~[N_{{(10^2)^p} 0}\cup (\bigcup_{q\in\omega}~
N_{{(10^2)^p} {(101)^q} 1^3})]\cup 
\{(10^2)^\infty\}\mbox{,}$$ 
which is a $\neg ~D_{2}(\boraone)$ set. Towards a 
contradiction, assume that $A^\infty$ is $D_{2}(\boraone)$: 
$${A^\infty=U_{1}\cap F=U\cup F_{2}}\mbox{,}$$ 
where the $U$'s are open and the $F$'s are closed. Let 
$O$ be a clopen set separating $\neg ~U_{1}$ from $F_{2}$ (see 22.C in 
[K1]). Then $A^\infty = (U\cap O)\cup (F\setminus O)$ 
would be in $\boraone\oplus\bormone$. If $(10^2)^\infty\in O$, then we would have 
$N_{(10^2)^{p_{0}}}\subseteq O$ for some integer $p_{0}$. But the 
sequence ${\big((10^2)^p}(1^20)^\infty\big)_{p\geq p_{0}}\subseteq 
O\setminus U$ and tends to $(10^2)^\infty$, which is absurd. If 
$(10^2)^\infty\notin O$, then we would have 
$N_{(10^2)^{q_{0}}}\subseteq\neg ~O$ for some integer $q_{0}$. But the 
sequence ${\big((10^2)^{q_{0}}} {(101)^q} 1^\infty\big)_{q\geq q_{0}}\subseteq 
F\setminus O$ and tends to ${(10^2)^{q_{0}}}(101)^\infty$, which is 
absurd.\bigskip

\noindent $\bullet$ For the class $D_2(\boraone)$, we can take $A := [A_{1}^{<\omega}]^{*}
\setminus [A_{0}^{<\omega}]^{*}$, where $A_{0}:= \{010,01^2\}$ and  
$$A_{1}:= \{010,01^2,0^2,0^3,10^2,1^20,10^3,1^20^2\}.$$ 
We have $A^\infty = A_1^\infty\setminus A_0^\infty$. Indeed, as $A\subseteq [A_{1}^{<\omega}]^{*}$, we have $A^\infty\subseteq A_{1}^\infty$. If $\alpha\in A_{0}^\infty$, then its decomposition into 
words of $A_{1}$ is unique and made of words in $A_{0}$. Thus $\alpha\notin A^\infty$ and 
$$A^\infty\subseteq A_1^\infty\setminus A_0^\infty .$$ 
Conversely, if $\alpha = {a_{0}}{a_{1}}
\ldots\in A_{1}^\infty\setminus A_{0}^\infty$, with $a_{i}\in A_{1}^{-}$, then 
there are two cases. Either there are infinitely many indexes $i$ 
(say $i_{0},i_{1},\ldots$) such that $a_{i}\notin A_{0}$. In this case, 
the words ${a_{0}}\ldots a_{i_{0}}$, ${a_{i_{0}+1}}
\ldots a_{i_{1}}$, $\ldots$, are in $A$ and $\alpha\in A^\infty$. 
Or there exists a maximal index $i$ such that $a_{i}\notin A_{0}$. In 
this case, ${a_{0}}\ldots {a_{i}} 0,~
10^2,~ 1^20\in A$, thus $\alpha\in A^\infty = A_1^\infty\setminus 
A_0^\infty$. Proposition 2 shows that 
$A\in D_2(\boraone)$. If $A^\infty = U\cup F$, with $U\in\boraone$ and 
$F\in\bormone$, then we have $U=\emptyset$ because $A_{1}^\infty$ is nowhere 
dense (every sequence in $A_{1}$ contains $0$, thus the sequences in $A_{1}^\infty$ 
have infinitely many $0$'s). Thus $A^\infty$ would be closed. But this 
contradicts the fact that $\big((01^2)^n0^\infty\big)_{n}\subseteq A^\infty$ 
and tends to $(01^2)^\infty\notin A^\infty$. Thus $A^\infty$ is a true 
$D_2(\boraone)$ set.\bigskip

\noindent $\bullet$ For the class $\check D_{3}(\boraone)$, we can take 
$A:= ([A_{2}^{<\omega}]^{*}\setminus [A_{1}^{<\omega}]^{*})\cup 
[A_{0}^{<\omega}]^{*}$, 
where $A_{0}:=\{ 0^2\}$, ${A_{1}:=\{ 0^2,01\}}$, and 
$A_{2}:=\{ 0^2,01,10,10^2\}$. We have $A^\infty = (A_{2}^\infty\setminus 
A_{1}^\infty )\cup A_{0}^\infty$. Indeed, as ${A\subseteq [A_{2}^{<\omega}]^{*}}$, 
we have ${A^\infty\subseteq A_{2}^\infty}$. If $\alpha\in A_{1}^\infty$, then 
its decomposition into words of $A_{2}^{-}$ is unique and made of 
words in $A_{1}$. If moreover $\alpha\notin A_{0}^\infty$, then it is clear 
that $\alpha\notin A^\infty$ and 
$${A^\infty \subseteq (A_{2}^\infty\setminus A_{1}^\infty )\cup A_{0}^\infty}.$$ 
Conversely, it is clear that 
$A_{0}^\infty\subseteq A^\infty$. If ${\alpha = {a_{0}} {a_{1}}
\ldots\in A_{2}^\infty\setminus A_{1}^\infty}$, then the argument above 
still works. We have to check that ${s:={a_{0}} \ldots  
a_{i_{0}}\notin [A_{1}^{<\omega}]^{*}}$. It is clear if ${a_{i_{0}}=10}$. 
Otherwise, $a_{i_{0}}=10^2$ and we argue by contradiction.

\vfill\eject

 The length of $s$ is even and the decomposition of $s$ into words of $A_{1}$ is unique. 
It finishes with $0^2$, and the even coordinates of the sequence $s$ are $0$. 
Therefore, $a_{i_{0}-1}\! =\! 0^2$ or $10$; we have the same thing with 
$a_{i_{0}-2}$, $a_{i_{0}-3},~\ldots$ Because of the parity, some $0$ 
remains at the beginning. But this is absurd. Now we have to check that 
${a_{0}}\ldots a_{i} 0\notin [A_{1}^{<\omega}]^{*}$. It is clear if $a_{i}=10^2$. 
Otherwise, $a_{i}=10$ and the argument above works.\bigskip 

 Finally, we have to check that if ${\gamma\in A_{1}^\infty}$, 
then ${\gamma -(0)\in A^\infty}$. There is a sequence $p_{0}, ~p_{1},\ldots$, 
finite or not, such that ${\gamma = ({0^{2p_{0}}})({{01}^{}})(  
{0^{2p_{1}}})( {{01}^{}})  {{\ldots}^{}} 0^\infty}$. Therefore 
$${\gamma -(0)= {({0^{2p_{0}}} 10)} {({0^{2p_{1}}} 10)}
{{\ldots}^{}}  (0^2)^\infty\in A^\infty}.$$

 If we set $U_{i}:=\neg ~A_{2-i}^\infty$, then we see that $A^\infty\in \check 
D_{3}(\boraone)$. If $\alpha$ finishes with $1^\infty$, then 
$\alpha\notin A_{2}^\infty$; thus $A_{2}^\infty$ is nowhere dense, 
just like $A^\infty$. Thus if $A^\infty = (U_{2}\setminus U_{1})\cup U_{0}$ 
with $U_{i}$ open, then $U_{0}=\emptyset$. By uniqueness of the decomposition 
of a sentence in $A_{i}^\infty$ into words of $A_{i+1}$, we see that 
$A_{i}^\infty$ is nowhere dense in $A_{i+1}^\infty$. So let 
$x_{\emptyset}\in A_{0}^\infty$, $(x_{n})\subseteq A_{1}^\infty\setminus 
A_{0}^\infty$ converging to $x_{\emptyset}$, and $(x_{n,m})_{m}\subseteq 
A_{2}^\infty\setminus A_{1}^\infty$ converging to $x_{n}$. Then 
$x_{n,m}\in U_{1}$, which is absurd. Thus 
$A^\infty\notin D_{3}(\boraone)$.\bigskip

\noindent $\bullet$ For the class $\check D_{2}(\boratwo)$, we can take 
${A:= \{s\!\in\! 2^{<\omega}~/~1^2\prec s~{\rm or}~s=(0)\}}$. We can write 
$${A^\infty = (\{0^\infty\}\cup\bigcup_{p} N_{0^p1^2})\cap (P_{f}
\cup\{\alpha\!\in\! 2^\omega /\forall n~\exists m\geq n~\alpha (m)\! =\! \alpha 
(m+1)\! =\! 1\})}.$$
Then $A^\infty\notin D_{2}(\boratwo)$, otherwise 
$A^\infty\cap N_{1^2}\in D_{2}(\boratwo)$ and would be a comeager 
subset of $N_{1^2}$. We could find $s\in 2^{<\omega}$ with even length 
such that $A^\infty\cap N_{{1^2} s}\in\bormtwo$. We define a continuous 
function $f:2^\omega\rightarrow 2^\omega$ by formulas 
$f(\alpha )(2n):=\alpha (n)$ if $n>\frac {\vert s\vert +1}{2}$, $({1^2} s)(2n)$ 
otherwise, and 
$f(\alpha )(2n+1):=0$ if $n>\frac {\vert s\vert }{2}$, $({1^2} s)(2n+1)$ otherwise. It 
reduces $P_{f}$ to $A^\infty\cap N_{{1^2}s}$, which is absurd.

\vfill\eject  

\bf\centerline{Summary of the complexity results in this paper:}\rm
$$\begin{tabular}{|c|c|c|c|c|}
\hline
& Baire category & complexity $\vert$ $\xi =1$ & $\xi = 2$ & $\xi\geq 3$\\
\hline
${\bf\Sigma}_{0}$ & nowhere dense & \multicolumn{3}{|c|}{${\Bormone\setminus\boraone}$}\\
\hline
${\bf\Pi}_{0}$ & co-nowhere dense & \multicolumn{3}{|c|}{${\Boraone\setminus\bormone}$}\\
\hline
${\bf{\it\Delta}}_{1}$ & co-nowhere dense & \multicolumn{3}{|c|}{${K_{{\it\sigma}}\setminus\bormtwo}$}\\
\hline
${\it\Sigma}_{\xi}$ & co-nowhere dense & \multicolumn{2}{|c|}{${\Ca\setminus\borel}$} & 
${\ca\setminus\borel}$\\
\hline
${\it\Pi}_{\xi}$ & co-nowhere dense & ${\Ca\setminus\bormtwo}$ & ${\Ca\setminus\borel}$ & 
${\ca\setminus\borel}$\\
\hline
${\it\Delta}$ & co-nowhere dense & \multicolumn{3}{|c|}{${\Ca\setminus\borel}$}\\
\hline
${\bf\Sigma}_{\xi}$ & co-nowhere dense & ${{\it\Delta}^1_{2}\setminus D_{2}(\boraone )}$ & 
${{\it\Sigma}^1_{2}\setminus \bormtwo}$ & ${{\bf\Sigma}^1_{2}\setminus D_{2}(\boraone )}$\\
\hline
${\bf\Pi}_{\xi}$ & co-nowhere dense & ${{\it\Delta}^1_{2}\setminus \bormtwo}$ & 
${{\it\Sigma}^1_{2}\setminus D_{2}(\boraone )}$ & ${{\bf\Sigma}^1_{2}\setminus D_{2}(\boraone )}$\\
\hline
${\bf\Delta}$ & co-nowhere dense & \multicolumn{3}{|c|}{${{\it\Sigma}^1_{2}\setminus D_{2}(\boraone )}$}\\
\hline
${\cal G}_{\xi}~(\xi\in\omega)$ &  & 
\begin{tabular}{c}
$\bormone\setminus\boraone$\\  nowhere dense
\end{tabular}
& ${\check D_{\omega}(\boraone )\setminus D_{\omega}(\boraone )}$ 
& ${\Ca\setminus D_{\omega}(\boraone )}$\\
\hline
${\cal F}$ &  & \multicolumn{3}{|c|}{${\Ca\setminus\bormtwo}$}\\
\hline
$\cal A_{\xi}$ & co-meager & 
\begin{tabular}{c}
$\check D_{2}(\boraone )\setminus D_{2}(\boraone )$\\  co-nowhere dense
\end{tabular}
& ${{\it\Sigma}^1_{2}\setminus D_{2}(\boraone )}$ 
& ${{\bf\Sigma}^1_{2}\setminus D_{2}(\boraone )}$\\
\hline
$\cal M_{\xi}$ & co-meager & 
\begin{tabular}{c}
${\it\Sigma}^1_{2}\setminus D_{2}(\boraone )$\\  co-nowhere dense
\end{tabular}
& ${{\it\Sigma}^1_{2}\setminus D_{2}(\boraone )}$ 
& ${{\bf\Sigma}^1_{2}\setminus D_{2}(\boraone )}$\\
\hline
${\cal B}$ &  & \multicolumn{3}{|c|}{${{\it\Sigma}^1_{2}}$}\\
\hline
${\cal A}$ & co-nowhere dense & \multicolumn{3}{|c|}{${{\it\Sigma}^1_{3}\setminus D_{2}(\boraone )}$}\\
\hline
\end{tabular}$$\smallskip

\noindent\bf {\Large 8 References.}\rm\bigskip

\noindent [F1]\ \ O. Finkel,~\it Borel hierarchy and omega context free languages,~\rm Theoret. Comput. Sci.~290, 3 (2003), 1385-1405

\noindent [F2]\ \ O. Finkel,~\it Topological properties of omega context free languages,~\rm Theoret. Comput. Sci.~262 (2001), 669-697

\noindent [K1]\ \ A. S. Kechris,~\it Classical Descriptive Set Theory,~\rm 
Springer-Verlag, 1995

\noindent [K2]\ \ A. S. Kechris,~\it On the concept of $\ca$-completeness,~\rm Proc. A.M.S.~
125, 6 (1997), 1811-1814

\noindent [L]\ \  A. Levy,~\it Basic set theory,~\rm Springer-Verlag, 1979

\noindent [Lo]\ \  M. Lothaire,~\it Algebraic combinatorics on words,~\rm Cambridge University Press, 2002

\noindent [Lou]\ \ A. Louveau,~\it A separation theorem for $\Ana$ sets,~\rm Trans. A. M. S.~260 (1980), 363-378

\noindent [M]\ \ Y. N. Moschovakis,~\it Descriptive set theory,~\rm North-Holland, 1980

\noindent [S]\ \ P. Simonnet,~\it Automates et th\'eorie descriptive,~\rm Ph. D. 
Thesis, Universit\'e Paris 7, 1992

\noindent [St1]\ \ L. Staiger,~\it $\omega$-languages,~\rm Handbook of Formal Languages, Vol 3, edited by G. Rozenberg and A. Salomaa, 
Springer-Verlag, 1997

\noindent [St2]\ \ L. Staiger,~\it On $\omega$-power languages,~\rm 
New Trends in Formal Languages, Control, Cooperation and Combinatorics, Lect. Notes in Comput. Sci. 1218 Springer-Verlag (1997), 377-393 

\vfill\eject

\end{document}